\documentclass[11pt]{article}
\usepackage[numbers]{natbib}
\usepackage{graphicx,amsmath,epsfig,amssymb,amsfonts,amsthm,enumerate,verbatim,todonotes,hyperref,mathtools, dsfont}

\usepackage{tikz,tikz-cd}
\usetikzlibrary{matrix,arrows}

\newtheorem{theorem}{Theorem}[section]
\newtheorem{corollary}[theorem]{Corollary}
\newtheorem{lemma}[theorem]{Lemma}

\newtheorem{question}[theorem]{Question}

\theoremstyle{definition}
\newtheorem{definition}[theorem]{Definition}

\newcommand{\restrict}{\mathord{\upharpoonright}}

\begin{document}

\title{Sets of real numbers closed under Turing equivalence:\\Applications to fields, orders and automorphisms}
\author{Iv\'an Ongay-Valverde\\\emph{Host Institution: Department of Mathematics}\\\emph{University of Wisconsin--Madison, Madison, WI, USA}\\\emph{Current Institution: Department of Mathematics and Statistics}\\\emph{York University, Toronto, ON, CA}\\\emph{Email: ivan.ongay.valverde@gmail.com} \\ ORCID: 0000-0002-7149-6061
}

\date{\begin{tabular}{rl}First Draft:&March 03, 2018\\Current Draft:&\today\end{tabular}}
\maketitle
\begin{abstract}
In the first half of this paper, we study the way that sets of real numbers closed under Turing equivalence sit inside the real line from the perspective of algebra, measure and order. Afterwards, we combine the results from our study of these sets as orders with a classical construction from Avraham to obtain a restriction about how non trivial automorphism of the Turing degrees (if they exist) interact with $1$-generic degrees.
\end{abstract}

{\small
\textbf{Keywords:} Computability Theory, Set Theory, Automorphism problem, Turing equivalence.

\textbf{MSC codes:} 03D28, 03D99, 03E75

\begin{center}
\textbf{Declarations}    
\end{center}

\textbf{Funding} Work done while being supported by CONACYT scholarship for Mexican students studying abroad and the Bank of Mexico through the FIDERH.

\textbf{Conflicts of interest/Competing interests} Not applicable

\textbf{Availability of data and material} Not applicable

\textbf{Code availability} Not applicable

\textbf{Authors' contributions} Not applicable}

\newpage

\section{Introduction}

Studying sets of reals has been one of the main objectives of set theory since it was conceived by Cantor (while studying derived sets). Impressively, as remarked by L\"owe \cite{BenediktTuringSetTh}, not much attention has been given to the study of sets of reals closed under Turing equivalence (except in a small list of articles mentioned in \cite{BenediktTuringSetTh}). Although it looks simple, this question can have many modifications that span a really rich area of study.

For example, focusing on algebra, we can ask which degrees are needed to form a field closed under Turing equivalence (or a vector space over $\mathds{Q}$). Furthermore, given $A\subseteq \mathds{R}$ we wonder how does the operator of `smallest field that contains $A$' and the operator `close $A$ under Turing equivalence' relate to each other. Do they conmute?

Focusing on measure we can ask which degrees are needed to form a measurable set that is closed under Turing equivalence. Finally, from the order perspective, we wonder which order types are attainable by sets closed under Turing equivalence.

All these questions are addressed in this work and, surprisingly, we discover that our work involving order types has an application to the automorphism problem of the Turing degrees. This work, as well as recent work in effective cardinal invariants, adds support to the old standing idea that studying computability related questions with set theoretic tools, and vice-versa, gives birth to new and interesting results in both areas.

The paper is organized as follows: the rest of this section will introduce some notation and give some observations that help explain why we choose to study certain questions in later sections. Sections \ref{Section algebra}, \ref{Section measure}, \ref{Section order} study sets of reals closed under Turing equivalence from the perspective of algebra, measure and order (respectively). In Section \ref{Decomposition Section} we prove that any Borel function is the countable union of monotone functions (Theorem \ref{Decomposition}). This theorem is a key component of our application to the automorphism problem.

Section \ref{Tower Section} centers on Theorem \ref{Automorphism thm}. This theorem is a restriction in the way that nontrivial automorphism of the Turing degrees interact with $1$-generic degrees. Finally, Section \ref{Modification Section} talks of possible ways to modify, or improve, Theorem \ref{Automorphism thm} and Section \ref{Question section} has open questions and conclusions.

I would like to warmly thanks Mariya Soskova for all his support and advice during this project and Joe Miller for all his useful comments and his permission to publish Lemma \ref{Decomposition counterexample}.

A version of this work appeared originally as a chapter of my doctoral dissertation \cite{IOVthesis}.

\section{Background and firsts definitions}

Given a real number $r$,  the set $\mathds{R}_{r}=\{s\in \mathds{R}: deg(r)=deg(s)\}$ is countable and dense. It is countable since we only have countably many Turing operators, it is dense since changing finite information does not change the Turing degree of a number. These two properties imply that $\mathds{R}_{r}$, as an order, is isomorphic to $\mathds{Q}$. 

On the other hand, given a noncomputable $r$, $\mathds{R}_{r}$ is not $\mathds{R}_{0}+r$ (since $2r\in \mathds{R}_{r}$ but $2r-r=r\notin \mathds{R}_{0}$) nor $r\cdot \mathds{R}_{0}$ (since $0\notin \mathds{R}_{r}$ and, less triavially, $r+\frac{1}{2}\in \mathds{R}_{r}$ but $1+\frac{1}{2r}\notin \mathds{R}_{0}$). So $\mathds{R}_{r}$ is not a coset of $\mathds{R}_{0}$ as a group (either with respect to addition or multiplication), it is not a subfield, ideal or a subgroup of $\mathds{R}$.

Given $\mathcal{S}$ a subset of the Turing degrees, our objective is to understand how does $\mathds{R}_{\mathcal{S}}=\bigcup_{deg(r)\in S}\mathds{R}_{r}$ sit inside $\mathds{R}$.

\subsection{Notation}

\begin{itemize}
    \item We will denote real numbers by lower case letters, $x, y, z, a, b, c$.
    \item We will denote Turing degrees by bold lower case letters, $\mathbf{x}, \mathbf{y}, \mathbf{z}, \mathbf{a}, \mathbf{b}, \mathbf{c} $.
    \item Upper case letters, $A,B, E, L$, will denote subsets of real numbers.
    \item $\varphi_{e}$ will denote the $e$-th Turing functional. If no oracle is expressed, we will assume that we are using a computable oracle.
    \item $\mathcal{D}$ is the set of all Turing degrees.
    \item Upper case calligraphic letters, $\mathcal{S}, \mathcal{A}, \mathcal{B}$, will denote subsets of the Turing degrees.
    \item $\mathds{R}_{x}$ ($\mathds{R}_{\mathbf{x}}$) are all the members of the real line, $\mathds{R}$, that are Turing equivalent to $x$ (resp. that have degree $\mathbf{x}$).
    \item Given a subset of reals $A$, $\mathds{R}_{A}=\bigcup_{x\in A}\mathds{R}_{x}$.
    \item Given a subset of Turing degrees $\mathcal{S}$, $\mathds{R}_{\mathcal{S}}=\bigcup_{\mathbf{x}\in \mathcal{S}}\mathds{R}_{\mathbf{x}}$.
\end{itemize}

\section{Algebra}\label{Section algebra}

The first question that we will work on is which subsets of the Turing degrees $\mathds{\mathcal{S}}$, make $\mathds{R}_{\mathcal{S}}$ a subfield of the real numbers.

From Theorem 6 in Rice \cite{rice1954recursive}, we know that $\mathds{R}_{0}$ is a real closed field. In other words, $\mathds{R}_{0}$ is an elementary subfield of the real numbers or, equivalently, a field that is linearly ordered and that has at least one root for all odd degree polynomials. Following the same proof as in \cite{rice1954recursive}, which boils down to showing that finding roots is a computable process, we can also show that $\mathds{R}_{\bigwedge x}$ is a real closed field, where $ \bigwedge x$ is the lower cone of $x$, i.e., all the reals that are computable from $x$.

\begin{corollary}
For any $x\in \mathds{R}$, the set $\mathds{R}_{\bigwedge x}$ is a real closed subfield of $\mathds{R}$.
\end{corollary}

Furthermore, we can be sure that all of these fields are non-isomorphic.

\begin{lemma}
Given $S, K\subseteq \mathds{R}$ real closed subfields of $\mathds{R}$ we have that $S$ can be embedded into $K$ if and only if $S\subseteq K$.
\end{lemma}

\begin{proof}
The ``if" side of the theorem is clearly true.

For the other implication, assume that $S$ can be embedded into $K$. Since both are subfields of the real numbers, we have that $\mathds{Q}\subseteq S\cap K$. Furthermore, any embedding between $S$ and $K$ will fix $\mathds{Q}$.

Now, since these are ordered fields, given $x\in \mathds{R}$, the real number is completely defined by it's order relation with respect to elements of $\mathds{Q}$. Therefore, the only embedding from $S$ to $K$ is inclusion.

\end{proof}

\begin{corollary}
The Turing degrees, with their order, can be embedded into the countable real subfields of $\mathds{R}$. 
\end{corollary}

Notice that $\mathds{R}_{\bigwedge x}$ is a subfield that contains $\mathds{R}_{x}$. Is it the smallest one? It turns out that it is.

\begin{lemma}\label{Bits and zeroes}
For any $x\in \mathds{R}$, and $y\in \mathds{R}_{\bigwedge x}$ there exists $z,w\in \mathds{R}_{x}$ such that $z+w\in \mathds{R}_{y}$.
\end{lemma}

\begin{proof}
Without loss of generality, we will assume that $x,y\in [0,1]$

Assume that we have sequences of digits $\langle a_{i} : i\in \omega \rangle $, $\langle b_{i} : i\in \omega \rangle $ that represents the decimal expansion of $x$ and $y$, respectively, and that are computable from $x$.

Define $s=\sum_{i=0}^{\infty}\frac{a_{i}}{10^{4i+1}}$ and $t=\sum_{i=0}^{\infty}\frac{b_{i}}{10^{4i+3}}$. Let $w=s+t$ and $z=-s$. Notice that you can easily distingish the decimal digits of $x$ and $y$ in $w$.

Clearly, $s$ and $t$ can be computed from $x$ and $y$ can compute $t$. Furthermore, we have that $t$ can compute $y$ and that $s$ and $w$ can compute $x$.

With all of this, we have that $w,z\in \mathds{R}_{x}$ and that $w+z=t\in \mathds{R}_{y}$.

\end{proof}

\begin{theorem}\label{Sum representation}
Given $y\in \mathds{R}_{\bigwedge x}$ we can express $y$ as a finite sum of elements of $\mathds{R}_{x}$.
\end{theorem}

\begin{proof}
Notice that, given $y\in [0,1]$ and $\langle b_{i} : i\in \omega \rangle $ that represents the decimal expansion of $y$ computable from $x$, we can write $y=a_{1}+a_{2}+a_{3}+a_{4}$ such that $a_{k}=\sum_{i=0}^{\infty}\frac{b_{4i+k}}{10^{4i+k}}$.

So, as we did in the proof of the above lemma, we can find $z_{k}, w_{k}\in \mathds{R}_{x}$ such that $a_{k}=z_{k}+w_{k}$. So, we have that $y=\sum_{k=1}^{4}z_{k}+w_{k}$.
\end{proof}

\begin{corollary}\label{Minimal field}
$\mathds{R}_{\bigwedge x}$ is the smallest field that contains $\mathds{R}_{x}$.
\end{corollary}

\begin{corollary}
If $\mathcal{C}$ is a set of Turing degrees that is linearly ordered under $\leq_{T}$ then $\mathds{R}_{\wedge \mathcal{C}}$ is the smallest field that contains $\mathds{R}_{\mathcal{C}}$. 
\end{corollary}

\begin{proof}
If $\mathcal{C}$ has a maximum element, then this follows from the result above. Otherwise, take $\langle a_{i} : i\in \alpha \rangle\subseteq \mathcal{C}$ a cofinal sequence, with $\alpha\in \{\omega, \omega_1\}$ and notice that $\bigcup_{i\in \alpha}\mathds{R}_{\wedge a_{i}}$ is a field.
\end{proof}

Now we can answer our first question:

\begin{theorem}
Given $\mathcal{S}$ a set of Turing degrees, we have that $\mathds{R}_{\mathcal{S}}$ is a field if and only if $\mathcal{S}$ is and ideal, i.e., $\mathcal{S}$ is closed under finite Turing joins (denoted by $\oplus$) and closed downward under Turing reduction.
\end{theorem}

\begin{proof}
In one direction, if $\mathds{R}_{\mathcal{S}}$ is a field, then given $x\in \mathds{R}_{\mathcal{S}}$ we have that $\mathds{R}_{x}\subseteq \mathds{R}_{\mathcal{S}}$. Therefore, by Corollary \ref{Minimal field}, we have that $\mathds{R}_{\wedge x}\subseteq \mathds{R}_{\mathcal{S}}$, so $\mathcal{S}$ is closed downwards. Now, if $x, y\in \mathds{R}_{\mathcal{S}}$ then, using an argument like the one in Lemma \ref{Bits and zeroes}, there are $s\in \mathds{R}_{x}$ and $t\in \mathds{R}_{y}$ such that $s+t=x\oplus y$.

For the other direction, assume that $\mathcal{S}$ is downward closed and closed under joins. First, notice that $0 , 1\in \mathds{R}_{\mathcal{S}}$ and that given $x\in \mathds{R}_{\mathcal{S}}$ we have that $\frac{1}{x}\in \mathds{R}_{\mathcal{S}}$. Finally, if $x, y\in \mathds{R}$ we have that $x+y, x\cdot y\leq_{T} x\oplus y$. So, given $x,y\in \mathds{R}_{\mathcal{S}}$, since $x\oplus y\in \mathcal{S}$, we can conclude that $x+y, x\cdot y\in \mathds{R}_{\wedge x\oplus y}\subseteq \mathds{R}_{\mathcal{S}}$.
\end{proof}

The following result shows that, when studying countable subfields of $\mathds{R}$ closed under Turing equivalence, it is enough to look at chains of Turing degrees.

\begin{corollary}
Given $F\subseteq \mathds{R}$ a countable subfield closed under Turing equivalence there is a set $\mathcal{C}$ of Turing degrees that is linearly order under $\leq_{T}$, such that $F=\mathds{R}_{\wedge \mathcal{C}}$.
\end{corollary}

\begin{proof}
Let $F=\{a_{i} : i\in \omega\}$ be a countable subfield closed under Turing equivalence and denote by $\textbf{a}_{i}$ the Turing degree corresponding to each element. Let $\mathcal{C}=\{\bigoplus_{i\in n} \textbf{a}_{i} : n\in \omega\}$. From the theorem above, we know that there is a set $\mathcal{S}$ of Turing degrees that is and ideal such that $F=\mathds{R}_{\mathcal{S}}$.

Notice that $\mathcal{C}\subseteq \mathcal{S}$ and that $\mathcal{S}\subseteq \wedge \mathcal{C}$. Therefore, $F=\mathds{R}_{\mathcal{S}}\subseteq \mathds{R}_{\wedge\mathcal{C}}$, $\mathds{R}_{\mathcal{C}}\subseteq \mathds{R}_{\mathcal{S}}=F$ and $\mathds{R}_{\wedge \mathcal{C}}\subseteq F$ (since $\mathds{R}_{\wedge \mathcal{C}}$ is the smallest field containing $\mathds{R}_{\mathcal{C}}$).
\end{proof}

In view of these result, we can conjecture that the smallest subfield containing $\mathds{R}_{x}\cup \mathds{R}_{y}$ is $\mathds{R}_{\wedge x\oplus y}$. Nevertheless, the following result hints to the contrary.

\begin{definition}
Given $A\subseteq \mathds{R}$ let $\langle A \rangle_{\mathds{Q}}$ be the smallest $\mathds{Q}$-vector space containing $A$. 
\end{definition}

\begin{corollary}
Given $x\in \mathds{R}$, $\langle \mathds{R}_{x}\rangle_{\mathds{Q}} =\mathds{R}_{\wedge x}$.
\end{corollary}

\begin{proof}
This follows from Theorem \ref{Sum representation}.
\end{proof}

\begin{theorem}
There are $x,y\in \mathds{R}$ such that $\langle \mathds{R}_{x}\cup \mathds{R}_{y}\rangle_{\mathds{Q}} \neq \mathds{R}_{\wedge x\oplus y} $.
\end{theorem}

\begin{proof}
First of all notice that, if $w\in \langle \mathds{R}_{x}\cup \mathds{R}_{y}\rangle_{\mathds{Q}}$ then there are $s\in \mathds{R}_{\wedge x}$ and $t\in \mathds{R}_{\wedge y}$ such that $w=s+t$. Then, if $w$ computes $s$ it follows that $w$ also computes $t$.

Using the classic result of embedability of upper semilattices as initial segments of the Turing degrees, by Lachlan and Lebeuf \cite{LachlanLebeuf}, there are degrees $\mathbf{x}, \mathbf{y}, \mathbf{w}$ such that $\mathbf{x}$ and $\mathbf{y}$ are distinct minimal degrees (the only degree below them is $\mathbf{0}$) and such that $\mathbf{x}<_{T} \mathbf{w} <_{T} \mathbf{x}\oplus \mathbf{y}$.

Notice that, given $w\in \mathbf{w}$, $w\neq s+t$ for any $s\in \mathds{R}_{\wedge x}$ and $t\in \mathds{R}_{\wedge y}$. To show this, suppose that $w=s+t$. Since $w$ computes $s$ we have that $w$ computes $t$. We know that $\mathbf{y}\not\leq_{T} \mathbf{w}$, therefore $t$ is computable. This means that $w=s+t$ is such that $w\in \mathds{R}_{x}$ which is impossible since $\mathbf{x}<_{T}\mathbf{w}$. This shows that $w\notin \langle \mathds{R}_{x}\cup \mathds{R}_{y}\rangle_{\mathds{Q}}$.

Furthermore, it also shows that $\mathds{R}_{\mathbf{w}}\cap \langle \mathds{R}_{x}\cup \mathds{R}_{y}\rangle_{\mathds{Q}}=\emptyset$ and that $\mathds{R}_{x\oplus y}\not\subseteq \langle \mathds{R}_{x}\cup \mathds{R}_{y}\rangle_{\mathds{Q}}$ (although, $x\oplus y\in \langle \mathds{R}_{x}\cup \mathds{R}_{y}\rangle_{\mathds{Q}}$). In particular, $\mathds{R}_{\langle \mathds{R}_{x}\cup \mathds{R}_{y}\rangle_{\mathds{Q}}}=\mathds{R}_{x\oplus y}\cup \mathds{R}_{x}\cup \mathds{R}_{y}\cup\mathds{R}_{0}\neq \mathds{R}_{\wedge x\oplus y}$, since $\mathds{R}_{\wedge x\oplus y}=\mathds{R}_{x\oplus y}\cup \mathds{R}_{w}\cup \mathds{R}_{x}\cup \mathds{R}_{y}\cup\mathds{R}_{0}$.

\end{proof}

\begin{corollary}
There is a collection of Turing degrees $\mathcal{S}$ such that $\langle \mathds{R}_{\mathcal{S}}\rangle_{\mathds{Q}}$ is not closed under Turing equivalence.
\end{corollary}

\begin{corollary}
There is a $Q$-vector space, $V$, inside $\mathds{R}$ such that $\mathds{R}_{V}$ is not a vector space.
\end{corollary}

Thanks to a suggestion from Josiah Jacobsen-Grocott, we can translate these results to the field case

\begin{theorem}\label{Josiah Subfield}
There are real numbers $x, y\in \mathds{R}$, $i\in \omega$  such that the minimal field containing $ \mathds{R}_x \cup \mathds{R}_y$ is not $\mathds{R}_{\wedge x\oplus y}$.
\end{theorem}

\begin{proof}
Let $x,y\in \mathds{R}$ be such that $ \mathbf{0}'= deg(x)'=deg(y)'=deg(x\oplus y)$. A construction of this pair of low degrees can be found making a small modification of the classic Sacks Splitting Theorem \cite{SacksSplit}.

Notice that the elements of the minimal field containing $\mathds{R}_x \cup \mathds{R}_y$ are of the form $\displaystyle \frac{\sum_{i=0}^{n}s_{i}t_{i}}{\sum_{j=0}^{m}a_{j}b_{j}}$ with $s_i, a_j\in  \mathds{R}_{\wedge x}$ and $t_i, b_j\in  \mathds{R}_{\wedge y}$.

Since $\mathbf{0}'=deg(x)'=deg(y)'$, we know that $\mathbf{0}'$ can compute a list $\{f_{e}\}_{e\in \omega}$ with all the $x$ computable reals, i.e., members of $\mathds{R}_{\wedge x} $ and a list $\{g_{e}\}_{e\in \omega}$ with all the $y$ computable reals. Therefore we can make a list of all $h_{e}=\displaystyle \frac{\sum_{i=0}^{n}f_{e_i}g_{h_i}}{\sum_{j=0}^{m}f_{k_j}g_{\ell_j}}$.

Using a diagonalization, there is a real $w$ computable from $\mathbf{0}'$ such that for all $e\in \omega$ $w\neq h_{e}$.

This means that $w\in \mathds{R}_{\wedge \mathbf{0}'}=\mathds{R}_{\wedge x\oplus y}$, but $w$ is not in the minimal field generated by $ \mathds{R}_x \cup \mathds{R}_y$.

\end{proof} 

 These few examples do not exhaust the questions related to $\mathds{R}_{\mathcal{S}}$ from an algebraic perspective. Nonetheless, we will leave some open questions and list them in the last section. Meanwhile, we will study these sets from a different perspective.

\section{Measure}\label{Section measure}

One of the firsts examples of pathological subsets of the real line appeared in the early 1900's when Guiseppe Vitali consruct a set that was not Lebesgue measurable. From that moment onward, finding strange subsets of the real line has been a usual part of the mathematical world. And, as we will see multiple times in this article, there are examples of sets of reals $A$ that have pathological properties but such that $\mathds{R}_{A}$ does not have the same properties. Because of this, we believe that `exploding' a set into a set that is closed under Turing equivalence is an operation that `tames' certain sets.

\begin{definition}
\begin{enumerate}
    \item A \emph{Vitali set}, $V\subseteq \mathds{R}$, is a set such that for any $r\in \mathds{R}$ there is $x\in V$ such that $x-r\in \mathds{Q}$ and given different $x,y\in V$, $x-y\notin \mathds{Q}$.
    \item A \emph{multiplicative Vitali set}, $M\subseteq \mathds{R}$, is a set such that for any $r\in \mathds{R}\setminus\{0\}$ there is $x\in M$ such that $\frac{x}{r}\in \mathds{Q}$ and given different $x,y\in M$, $\frac{x}{y}\notin \mathds{Q}$.
\end{enumerate}
\end{definition}

\begin{theorem}[Folklore]
Vitali sets and multiplicative Vitali sets are not measurable.
\end{theorem}

\begin{theorem}
Given $V$ a Vitali set, $\mathds{R}_{V}=\mathds{R}$. Analogously, given $M$ a  multiplicative Vitali set, $\mathds{R}_{M}=\mathds{R}$.
\end{theorem}

\begin{proof}
We will use that addition and multiplication by a computable number are computable operations.

In the additive case, we know that $\mathds{R}=V+\mathds{Q}=\{r+q: q\in \mathds{Q}, r\in V\}$ and that $r+q$ is Turinq equivalent to $r$. Therefore, $\mathds{R}_{V}=\mathds{R}$.

We can do an analogue argument for the multiplicative case.
\end{proof}

This shows that there are non measurable sets that become measurable once they are closed under Turing equivalence.

In the literature we can find a direct generalization of Vitali sets focussing on groups different than $\mathds{Q}$. These examples will behave differently.

\begin{definition}
\begin{enumerate}
    \item Given a group $G$ and an action $a$ of $G$ over $\mathds{R}$, we say that $S\subseteq \mathds{R}$ is an \emph{$a$-selecto}r (or \emph{$G$-selector} when the action is understood) if and only if for every $x\in \mathds{R}$ there is $v\in S$ and $g\in G$ such that $g\cdot x=v$ and for all $x,y\in S$, $x\neq y$, there is no $g\in G$ such that $g\cdot x=y$.
    \item We say that an action $a$ of $G$ over $\mathds{R}$ is \emph{paradoxical} if and only if every selector of $a$ is non measurable.
    \item If $a$, an action of $G$ over $\mathds{R}$, is paradoxical we call the selectors \emph{$a$-Vitali sets} and, in case that the action is clear, we call them \emph{$G$-Vitali sets}.
\end{enumerate}
\end{definition}

\begin{lemma}[Folklore]
Given a measurable set $A$, $\displaystyle \lim_{r\rightarrow 0} \lambda(A\cap (A+r))=\lambda(A)$ where $\lambda$ is the Lebesgue measure (in any dimension) and $A+r=\{a+r: a\in A\}$.
\end{lemma}

\begin{lemma}
For any computable group $G$, i.e. a countable group with a computable bijection with $\omega$ such that the group operation is computable, and a computable action $a$ of $G$ over $\mathds{R}$ we have that for every $a$-selector set $A$,  $\mathds{R}_{A}=\mathds{R}$.
\end{lemma}

\begin{theorem}
There is a group $G$ such that for all its $G$-Vitali sets, $A$, $\mathds{R}_{A}$ is nonmeasurable.
\end{theorem}

\begin{proof}

Let $r\in \mathds{R}$ be noncomputable and let $G=r\cdot \mathds{Q}=\{rq: q\in \mathds{Q}\}$. This set is a group under addition and it has a cannonical action over $\mathds{R}$ (again, addition). Since $G$ is isomorphic to $\mathds{Q}$ as a group, the same proof that shows that any $\mathds{Q}$-selector is nonmeasurable shows that any $G$-selector is nonmeasurable.

Now, let $A$ be a $G$-Vitali set. We want to show that $\mathds{R}_{A}$ is nonmeasurable. In order to do this, we will show that $\mathds{R}_A\times \mathds{R}_A\subseteq \mathds{R}^{2}$ is nonmeasurable.

Notice that, since $A+G=\mathds{R}$ we also have that $\mathds{R}_{A}+G=\mathds{R}$, then $\mathds{R}_{A}\times \mathds{R}_{A}+G\times G=\mathds{R}^{2}$. This means that, $\mathds{R}_{A}\times \mathds{R}_{A}$ is not of measure zero (since $G$ is countable).

To finish the proof, we would like to show that $\mathds{R}_{A}\cap\left(\mathds{R}_{A}+qr\right)=\emptyset$. Nevertheless, this is not true. Given $x\in \mathds{R}_{A}$ and $q\in \mathds{Q}$ if $x+qr=y\in \mathds{R}_{A}$ then  $x\oplus y$ computes $r$. We will have to modify $\mathds{R}$ by a null set to finish our proof. 

Let $B=\mathds{R}_{A}\times \mathds{R}_{A}\setminus N$ where \[N=\{(x,y) : r\leq_{T} x\oplus y\}. \]
The classical proof of Sacks that the upper cone of non-computable reals have measure zero in $\mathds{R}$, in \cite{SacksDeg}, also shows that $N$ is a null set of $\mathds{R}^{2}$ (for non computable $r$). It is enough to show that $B$ is nonmeasurable to finish the proof.

Notice that, given $h\in G\times G\setminus \{(0,0)\}$, $(B+h)\cap B=\emptyset$. Therefore, for any sequence from $G\times G$, $h_{i}$, that converges to $(0,0)$ we have that $\displaystyle\lim_{i\rightarrow \infty}\lambda((B+h_{i})\cap B)=0$. If $B$ was a measurable set, it would have to have measure $0$, but we know that $B$ is not a null set. Therefore, $B$ is not measurable.

\end{proof}

\begin{corollary}\label{Cantor set trick}
There is a measurable set $L$ such that $\mathds{R}_{L}$ is nonmeasurable.
\end{corollary}

\begin{proof}
Let $C$ be the ternary Cantor set. Take a function, $f$, that bijects $2^{\omega}$ with the usual Cantor set inside the reals, $C$, and maintains degree, in other words, given $a\in 2^{\omega}$ then $deg(a)=deg(f(a))$.

This functions shows that $C$ has a representative of every degree in $\mathds{R}$.
Using the set $A$ from the theorem above, let $L=\mathds{R}_{A}\cap C$.

Since $C$ has measure zero we know that $L$ is a measurable set of measure zero. Now, since for all real number there is a Turing equivalent real in $C$ we have that $\mathds{R}_{L}=\mathds{R}_{A}$.
\end{proof}

Although it would have been interesting to find a positive measure measurable set whose closure under Turing equivalence is not measurable, we found out that this is not possible thanks to a suggestion from James Hanson.

\begin{theorem}\label{Lebesgue positive}
Given $A\subseteq \mathds{R}$ such that $\lambda(A)>0$, we have that $\mathds{R}_{A}$ is measurable.
\end{theorem}

\begin{proof}
Remember that the Lebesgue Density Theorem states that given a measurable set $A$, for almost all points of $A$, $\displaystyle\lim_{\epsilon \rightarrow 0} \frac{\lambda (A\cap (x-\epsilon, x+\epsilon))}{\lambda((x-\epsilon, x+\epsilon))}$ is either $0$ or $1$. Furthermore, if it it the case that for almost all points the Lebesgue density is $0$, then $A$ is a measure zero set.

We will show that if $\lambda (A)>0$ then $\mathds{R}\setminus(\mathds{Q}+A)$ has measure $0$, i.e., $\mathds{Q}+A$ has full measure. Since $\mathds{Q}+A\subseteq \mathds{R}_{A}$, we have that $\mathds{R}_{A}$ has also full measure, hence, it is measurable.

Now, to show that $\mathds{Q}+A$ has full measure we will show that $\lambda((\mathds{Q}+A)\cap B)> 0$ for all $B$ such that $\lambda(B)>0$.

Let $B$ be such that $\lambda(B)>0$. Using the Lebesgue Density Theorem, find $a\in A$, $b\in B$ and $\epsilon>0$ such that $\displaystyle \lambda (A\cap (a-\frac{\epsilon}{2}, a+\frac{\epsilon}{2}))>\frac{3\epsilon}{4}$ and $\displaystyle \lambda (B\cap (b-\frac{\epsilon}{2}, b+\frac{\epsilon}{2}))>\frac{3\epsilon}{4}$.

Now, let $r\in \mathds{Q}$ such that $\displaystyle\lambda((a+r-\frac{\epsilon}{2}, a+r+\frac{\epsilon}{2})\cap (b-\frac{\epsilon}{2}, b+\frac{\epsilon}{2}))>\frac{3\epsilon}{4}$. Since each of the intervals lost at most $\displaystyle \frac{\epsilon}{4}$ lenght when intersected, we have that
\[\displaystyle\lambda((A+r)\cap(a+r-\frac{\epsilon}{2}, a+r+\frac{\epsilon}{2})\cap (b-\frac{\epsilon}{2}, b+\frac{\epsilon}{2}))>\frac{2\epsilon}{4}\]
and
\[\displaystyle\lambda(B \cap (a+r-\frac{\epsilon}{2}, a+r+\frac{\epsilon}{2})\cap (b-\frac{\epsilon}{2}, b+\frac{\epsilon}{2}))>\frac{2\epsilon}{4}.\]
Therefore,
\[\displaystyle\lambda((A+r)\cap B\cap (a+r-\frac{\epsilon}{2}, a+r+\frac{\epsilon}{2})\cap (b-\frac{\epsilon}{2}, b+\frac{\epsilon}{2}))>\frac{\epsilon}{4}.\]

This shows that $\lambda((A+r)\cap B)>0$ which implies that $\lambda((A+\mathds{Q})\cap B)>0$. 
\end{proof}

We think that it is important to investigate which non-measurable sets become measurable once they are closed under Turing equivalence.

Because it is well known that countable sets have measure zero and that, depending on set theory axioms, uncountable sets of size strictly less than the continuum are either measure zero or non measurable, while working on this topic we should focus our attention exclusively on sets of size $\mathfrak{c}$ (continuum). Nevertheless, in this work we will not investigate other examples related to measure.

\section{Order types}\label{Section order}

A different question about $\mathds{R}_{\mathcal{S}}$ is whether it can be isomorphic to any order type that sits inside the reals.

A really quick answer for this is ``no''. To see this, we need a definition:

\begin{definition}
Given $A\subseteq \mathds{R}$ and a cardinal $\kappa\leq 2^{\aleph_{0}}=\mathfrak{c}$ we say that $A$ is \emph{$\kappa$-dense} if and only if for any nonempty interval $I$ we have that $|I\cap A|=\kappa$. 
\end{definition}

Necessarily, $\mathds{R}_{\mathcal{S}}$ is going to be a $|\mathcal{S}|\cdot \aleph_{0}$-dense subset of reals (because each $\mathds{R}_{x}$ is $\aleph_{0}$-dense). Therefore, any order that is not $\kappa$-dense is not going to have an isomorphic copy as $\mathds{R}_{\mathcal{S}}$.  

By a classical result of Cantor, we know that every $\aleph_{0}$-dense subset of reals is order isomorphic to $\mathds{Q}$. Now, if there is $\aleph_{0}<\kappa<\mathfrak{c}$, then the study of $\kappa$-dense subsets of reals is completely independent of ZFC (see Baumgartner \cite{baumgartner1973alphall}). So we will only focus on $\mathfrak{c}$-dense subsets of the reals (or subsets such that $|\mathcal{S}|=\mathfrak{c}$).

Even in this situation, we know that we cannot get all possible order types, for example $\mathds{R}\setminus \{0\}$ will not be isomorphic to any $\mathds{R}_{\mathcal{S}}$ (any $\mathds{R}_{\mathcal{S}}$ that contains an interval will be $\mathds{R}$).

Trying to characterize which order types can be attained is still an open question. Towards giving more insight into the problem, we will focus our attention to specific $\mathfrak{c}$-dense subsets of reals.

\subsection{Luzin and Sierpinski sets}

\begin{definition}
We say that $A\subseteq \mathds{R}$ is a \emph{Luzin set} if and only if $A$ is uncountable and the intersection of $A$ with any meager set (equivalently, no-where-dense set) is countable.

$A$ is a \emph{Sierpinski set} if, in the above definition, we replace meager with null (or measure zero) set.
\end{definition}

The following lemma also appears in L\"owe \cite{BenediktTuringSetTh} with a more complicated proof, since they work in $2^{\omega}$ instead of the real line. Although we are sure that the results are the same, this proof shows how different it is to work with sets of real numbers depending on if you see them as contained in $2^{\omega}$, $\omega^{\omega}$ or the real line.

\begin{lemma}
There is no $\mathcal{S}$ such that $\mathds{R}_{\mathcal{S}}$ is a Luzin or Sierpinski set.
\end{lemma}

\begin{proof}
Let $C$ be the ternary Cantor set. As we remark in Corollary \ref{Cantor set trick}, $C$ has a representative of all Turing degrees in $\mathds{R}$. 

Notice that, $|\mathds{R}_{\mathcal{S}}\cap C|=|\mathcal{S}|$. Then, if $\mathcal{S}$ is countable, $\mathds{R}_{\mathcal{S}}$ is also countable, so it is not Luzin; and if $\mathcal{S}$ is uncountable, $\mathds{R}_{\mathcal{S}}$ intersects a meager set, $C$, in uncountably many points, so it is not Luzin. 

An analogous reasoning shows that it cannot be Sierpinski using that $C$ is a null set.
\end{proof}

We believe that it is important to remark that the above result is not as trivial as it seems; notice that:

\begin{lemma}\label{L+Q}
Given a Luzin set $L$ (resp. a Sierpinski set), the set $L+\mathds{Q}$ is also a Luzin set (resp. a Sierpinski set).
\end{lemma}

\begin{proof}
Let $L$ be an uncountable set and assume that $L+\mathds{Q}=\{l+r: l\in L, r\in \mathds{Q}\}$ is not Luzin. Therefore, there is a meager set $M$ such that $M\cap (L+\mathds{Q})$ is uncountable.

Define $L+r$ as $\{l+r: l\in L\}$. Since $\mathds{Q}$ is countable, there should be $r\in \mathds{Q}$ such that $(L+r)\cap M$ is uncountable.

Nevertheless,  $(M-r)\cap L$ is uncountable since $|(M-r)\cap L|=|(L+r)\cap M|$ and $M-r$ is also meager. This means that $L$ is not Luzin.

In the case for Sierpinski, the argument is the same but replacing $M$ with a null set.
\end{proof}

Taking into account that every uncountable subset of a Luzin (resp. Sierpinski) set is also Luzin (resp. Sierpinski) set. The above result can be improved to:

\begin{corollary}
Given a set $L\subseteq \mathds{R}$, we have that $L$ is a Luzin set (resp. a Sierpinski set) if and only if the  $L+\mathds{Q}$ is a Luzin set (resp. a Sierpinski set).
\end{corollary}

Analyzing the proof of the lemma above we can notice that the key point is that given $M$ meager, $M-r$ is also meager. If we think of translation by $r$ as a function, say $f$, we have that for all meager sets $M$, $f^{-1}(M)$ is also meager.

This observation lets us expand our lemma:

\begin{corollary}\label{Preserving Luzin}
Given a countable set $\mathcal{F}$ of functions whose inverses preserve category (resp. measure) closed under composition and a Luzin set (resp. Sierpinski set) $L$, we have that the closure of $L$ under $\mathcal{F}$ is also a Luzin set (resp. Sierpinski set).
\end{corollary}

\begin{proof}
We can follow the same proof as in Lemma \ref{L+Q} using the comments made in the observation above.
\end{proof}

\begin{corollary}\label{No Measure generation}
It is consistent with ZFC that there is no countable set $\mathcal{F}$ of functions whose inverses preserve category (resp. measure) closed under composition such that for all reals, $r$, we have that $\mathds{R}_{r}$ is the closure  under $\mathcal{F}$ of the set $\{r\}$.
\end{corollary}

\begin{proof}
Given a countable collection of functions, $\mathcal{F}$, call $\mathcal{F}_{r}$ the closure of $\{r\}$ under $\mathcal{F}$. Notice that, since $\mathcal{F}$ is countable, $\mathcal{F}_{r}$ is also countable. Given $S\subseteq \mathds{R}$, we define analogously $\mathcal{F}_{S}$.

It is consistent with ZFC, regardless of the cardinal arithmetic, that there exists a Luzin set $L$ (of any desired size).

If $\mathcal{F}$ is a countable set of functions whose inverses preserve category and that is closed under composition. Then, by the above corollary, we have that $\mathcal{F}_{L}$ is Luzin.

We can show that it is not the case that for all $r\in \mathds{R}$ we have that $ \mathds{R}_{r}\subseteq \mathcal{F}_{r}$: assume otherwise, then we will have that $\mathds{R}_{L}\subseteq \mathcal{F}_{L}$. Nevertheless, $\mathcal{F}_{L}$ is Luzin and $\mathds{R}_{L}$ is not, but this is a contradiction. Therefore, there exist uncountably many $r\in L$ such that $\mathcal{F}_{r}\neq \mathds{R}_{r}$.

The proof for measure preserving functions is analogous, but it assumes the existence of a Sierpinski set instead of a Luzin set.

\end{proof}

\begin{corollary}
It is consistent with ZFC that there is no countable set $\mathcal{F}$ of functions whose inverses preserve category (resp. measure) closed under composition such that for all reals, $r$, we have that $\mathds{R}_{\wedge r}$ is the closure  under $\mathcal{F}$ of the set $\{r\}$.
\end{corollary}

This last corollary can also be shown directly (in ZFC) by constructing a computable function that does not preserve category or measure. For the measure case, a candidate is a computable isomorphism between a measure zero Cantor set and a positive measure Cantor set.

\subsection{Entangled sets}

\begin{definition}
Given a collection of $2$-tuples $\{(x_{\alpha}, y_{\alpha}): \alpha\in I\}$ we say that these tuples are \emph{disjoint} if and only if $|\{x_{\alpha}, y_{\alpha}\}|=2$ for all $\alpha$ and $\{x_{\alpha}, y_{\alpha}\}\cap \{x_{\beta}, y_{\beta}\}=\emptyset$ for all $\alpha\neq \beta$.

Given $\aleph_{0}<\kappa\leq \mathfrak{c}$, we say that $A\subseteq \mathds{R}$ is \emph{$\kappa$-$2$-entangled} if and only if $|A|\geq\kappa$ and for every collection of size $\kappa$ of disjoint $2$-tuples of $A$ there are $(x_{1}, y_{1})$, $(x_{2}, y_{2})$, $(w_{1}, z_{1})$, $(w_{2}, z_{2})$ such that $x_{1}<x_{2}$ and $y_{1}<y_{2}$, $w_{1}<w_{2}$ but $z_{1}>z_{2}$.

If a set is $\aleph_{1}$-$2$-entangled we just say that it is \emph{$2$-entangled}.

\end{definition}

Notice that, any one-to-one function between two $\kappa$ size disjoint subsets of a $\kappa$-$2$-entangled set can be seen as a collection of size $\kappa$ of disjoint $2$-tuples. Then, by the defining property of $\kappa$-$2$-entangled sets, this one-to-one function cannot be monotone increasing nor monotone decreasing. This implies that, given any two disjoint subsets of size $\kappa$ inside of a $\kappa$-$2$-entangled, the subsets are not order isomorphic. 

\begin{lemma}
There is no uncountable $\mathcal{S}$ such that $\mathds{R}_{\mathcal{S}}$ is $|\mathcal{S}|$-$2$-entangled.
\end{lemma}

\begin{proof}
Notice that, for any $\mathcal{S}$, $\mathds{R}_{\mathcal{S}}+2=\mathds{R}_{\mathcal{S}}$, therefore, $\left([0,1]\cap \mathds{R}_{\mathcal{S}}\right)+2\subseteq \mathds{R}_{\mathcal{S}}$ and $\left([0,1]\cap \mathds{R}_{\mathcal{S}}\right)+2\subseteq [2,3]$. Which means that $\left(\left([0,1]\cap \mathds{R}_{\mathcal{S}}\right)+2 \right)\cap \left([0,1]\cap \mathds{R}_{\mathcal{S}}\right)=\emptyset$.

If $\mathcal{S}$ is uncountable, we have that $[0,1]\cap \mathds{R}_{\mathcal{S}}$ is uncountable. So the function $x\mapsto x+2$ is a one-to-one order preserving function\footnote{An order preserving function is a function such that $x<y$ then $f(x)<f(y)$. On the other hand, an order reversing function is such that if $x<y$ then $f(y)<f(x)$.} that shows that $\left([0,1]\cap \mathds{R}_{\mathcal{S}}\right)+2$ and $[0,1]\cap \mathds{R}_{\mathcal{S}}$ are order isomorphic.

This shows that $\mathds{R}_{\mathcal{S}}$ is not $|\mathcal{S}|$-$2$-entangled.
\end{proof}

Upon closer inspection of the proof above, we can see that what makes $\mathds{R}_{\mathcal{S}}$ not $|\mathcal{S}|$-$2$-entangled is that there are order preserving functions that preserve the degree. Making a slight modification to the definition of an entangled set, we can create an $\mathcal{S}$ such that $\mathds{R}_{\mathcal{S}}$ is almost $|\mathcal{S}|$-$2$-entangled.

\begin{definition}
Given $\aleph_{0}<\kappa\leq \mathfrak{c}$, we say that $A\subseteq \mathds{R}$ is \emph{layered $\kappa$-$2$-entangled} (or \emph{$\sigma$-$2$-entangled} if the $\kappa$ is understood) if and only if $|A|=\kappa$ and there is a function $ht:A\rightarrow \kappa$ countable-to-one such that for every collection of size $\kappa$ of disjoint $2$-tuples of $A$ such that the entries of a given pair have different values of $ht$ then there are $(x_{1}, y_{1})$, $(x_{2}, y_{2})$, $(w_{1}, z_{1})$, $(w_{2}, z_{2})$ such that $x_{1}<x_{2}$ and $y_{1}<y_{2}$, $w_{1}<w_{2}$ but $z_{1}>z_{2}$.
\end{definition}

Notice that any layered $\kappa$-$2$-entangled set contains a $\kappa$-$2$-entangled: every collection of $\kappa$ many elements that have different values of $ht$ is $\kappa$-$2$-entangled.

\begin{lemma}
If $A$ is $\kappa$-$2$-entangled and $|A|=\kappa$ then $A+\mathds{Q}$ is layered $\kappa$-$2$-entangled.
\end{lemma}

\begin{proof}
Let $A=\{x_{\alpha}: \alpha<\kappa\}$, given $s\in A+\mathds{Q}$ let $ht(s)=\min\{\alpha:\exists r\in \mathds{Q}\ (s=x_{\alpha}+r)\}$. Since $\mathds{Q}$ is countable, this function is countable-to-one. This function makes $A+\mathds{Q}$ layered $\kappa$-$2$-entangled.

To show this, assume that $A+\mathds{Q}$ is not layered $\kappa$-$2$-entangled. Then, there is a collection of size $\kappa$ of disjoint $2$-tuples of $A+\mathds{Q}$ such that the entries of a given pair have different values of $ht$ such that, as a function, call it $f$, from a subset of $A+\mathds{Q}$ to $A+\mathds{Q}$ it is order preserving or order reversing.

Without loss of generality, let's say that it is order preserving. Now, given any element $a\in dom(f)$ there is $r_{a}\in\mathds{Q}$ such that $a=x_{ht(a)}+r_{a}$. We know that $|dom(f)|=\kappa$, so, shrinking $f$, we can find a single $r$ such that for all $a\in dom(f)$, $a=x_{ht(a)}+r$. Since $f$ is one-to-one, we can do the same with $range(f)$, i.e., we can shrink $f$ in such a way that there is $q\in\mathds{Q}$ such that for all $b\in range(f) $ we have  $b=x_{ht(b)}+q$.

Finally, notice that the function $f(x+r)-q$ is an order preserving one-to-one function from $A$ to itself with domain of size $\kappa$. This contradicts the fact that $A$ is $\kappa$-$2$-entangled. 
\end{proof}

\begin{corollary}\label{From ent to lay ent}
Given a countable set $\mathcal{F}$ of strictly monotone functions closed under composition and a $\kappa$-$2$-entangled set $A$ of size $\kappa$ we have that the closure of $A$ under $\mathcal{F}$ is layered $\kappa$-$2$-entangled.
\end{corollary}

\begin{proof}
The proof is the same as above as long as we replace addition by a rational with other one-to-one monotone function.
\end{proof}

From now on, we will fix a bijection between the Turing degrees and $\mathfrak{c}$. While working with collections of Turing degrees our function $ht$ will be the composition of the function that maps a point to its degree followed by the bijection fixed above. In case that the size of our set is not $\mathfrak{c}$, we create a bijection with the corresponding cardinal. Nevertheless, to simplify notation, we will refer to $ht$ as the degree of a point.

Given an entangled set $A$, is it true that $\bigcup_{a\in A}\mathds{R}_{a}=\mathds{R}_{A}$ is layered entangled? At this moment we can show this to be true using the, really strong, Proper Forcing Axiom (PFA, see Chapter 5 section 7 of \cite{kunen2014set}). Later, using Theorem \ref{Decomposition}, we will see that PFA was not necessary for this result. Nevertheless, this first approach helps to show which steps need to be followed.

\begin{corollary}
$PFA$ implies that for every $\kappa$-$2$-entangled set, $A$, we have that $\mathds{R}_{A}$ is layered $\kappa$-$2$-entangled.
\end{corollary}

\begin{proof}
$PFA$ implies that any countable-to-one function between subsets of reals is the union of countably many monotone functions. Now, given a computable function $\varphi_{e}:\omega\rightarrow \omega$ we can associate it to $f_{e}:\mathds{R}\rightarrow \mathds{R}$ where $f_{e}(x)=\varphi^{x}_{e}$ (whenever $\varphi^{x}_{e}$ is total). Now, let $D_{e}=\{x: x\equiv_{T}f_{e}(x)\}$. We have that $f_{e}\restrict_{D_{e}}$ is countable-to-one.

Notice that, using the notation of \ref{No Measure generation}, $\mathcal{F}=\{f_{e}\restrict D_{e}: e\in \omega\}$ is such that for all $x\in \mathds{R}$, $\mathcal{F}_{x}=\mathds{R}_x$. Using $PFA$, we can change $\mathcal{F}$ for a countable collection of monotone functions. Combining this fact with Corollary \ref{From ent to lay ent} we are done.

\end{proof}

Later we will show in Theorem \ref{Turing is layered entangled} that $ZFC$ implies that given a $2$-entangled set $E$, $\mathds{R}_{E}$ is layered-$2$-entangled set. Nevertheless, it is important to remark that we can get certain control over a layered entangled set with the following techniques:

Recall that, given $\mathcal{S}$ a subset of the Turing degrees we define $\mathds{R}_{\mathcal{S}}=\bigcup_{deg(a)\in \mathcal{S}}\mathds{R}_{a}$. We will write Turing degrees with bold case to differentiate them from real numbers.

\begin{theorem}\label{continuum entangled}
Given continuum many Turing degrees, $\mathcal{A}$, there is a continuum size collection of degrees $\mathcal{S}\subseteq \mathcal{A}$ such that $\mathds{R}_{\mathcal{S}}$ is layered $\mathfrak{c}$-$2$-entangled.
\end{theorem}

\begin{proof}
Let $\{f_{\alpha}: \alpha<\mathfrak{c}\}$ be an enumeration of all continuous functions from $G_{\delta}$ subsets of $\mathds{R}$ to $\mathds{R}$ and let $\mathcal{A}=\{\mathbf{a}_{\alpha} :\alpha<\mathfrak{c}\}$ be an enumeration of continuum many Turing degrees. Given a Turing degree $\mathbf{a}$, we will denote by $\mathds{R}_{\mathbf{a}}$ all reals that have Turing degree $\mathbf{a}$.

We will construct the set $\mathcal{S}\subseteq \mathcal{A}$ by recursion.

Assume that we already have $\mathcal{S}_{\alpha}=\{\mathbf{a}_{\gamma_{\xi}}: \xi<\alpha\}$. Let \[\gamma_{\alpha}=\min\{\beta : \mathds{R}_{\mathbf{a}_{\beta}}\cap(\bigcup_{\xi<\alpha}\mathds{R}_{\mathbf{a}_{\gamma_{\xi}}}\cup B_{\alpha})=\emptyset\}\]
where \[B_{\alpha}=\{f_{\xi}(\overline{v}): \xi<\alpha, \overline{v}\in \bigcup_{\xi<\alpha}\mathds{R}_{\mathbf{a}_{\gamma_{\xi}}}\}.\]

We know that $\mathds{R}\setminus (\bigcup_{\xi<\alpha}\mathds{R}_{\mathbf{a}_{\gamma_{\xi}}}\cup B_{\alpha})$ is of size continuum because \[\left|\bigcup_{\xi<\alpha}\mathds{R}_{\mathbf{a}_{\gamma_{\xi}}}\cup B_{\alpha}\right|=|\alpha|<\mathfrak{c}.\] This means that there are at most $|\alpha|$  ordinals $\beta$ such that $\mathds{R}_{\mathbf{a}_{\beta}}\cap(\bigcup_{\xi<\alpha}\mathds{R}_{\mathbf{a}_{\gamma{\xi}}}\cup B_{\alpha})\neq\emptyset$ (this is because the Turing equivalent classes are disjoint between them).

The set $\mathcal{S}=\mathcal{S}_{\mathfrak{c}}=\{\mathbf{a}_{\beta_{\alpha}}: \alpha<\mathfrak{c}\}$ is the one that we are looking for. To show that $\mathds{R}_{\mathcal{S}}$ has the property that for every collection of size continuum of disjoint $2$-tuples of $\mathds{R}_{\mathcal{S}}$ such that each entry of each pair has different degree then there are   $(x_{1}, y_{1})$, $(x_{2}, y_{2})$, $(w_{1}, z_{1})$, $(w_{2}, z_{2})$ such that $x_{1}<x_{2}$ and $y_{1}<y_{2}$, $w_{1}<w_{2}$ but $z_{1}>z_{2}$, we will follow a technique used by Todorcevic in \cite{todorcevic1989partition}. We reproduce the argument here for the convenience of the reader.

Let $\{(x_{\alpha_{\xi}}, x_{\beta_{\xi}}): \xi<\mathfrak{c}\}\subseteq \mathds{R}_{\mathcal{S}}^{2}$ be a collection of size continuum of disjoint $2$-tuples of $\mathds{R}_{\mathcal{S}}$ such that each entry of each pair has different degree.

Let \[K=\{x_{\alpha_{\xi}}: \alpha_{\xi}< \beta_{\xi}<\mathfrak{c}\}.\] We can assume that this set is of size continuum. If not, we can run the argument interchanging the roles of $\alpha_{\xi}$ and $\beta_{\xi}$.

Now, we can define the function $g:K\rightarrow \mathds{R}$ such that \[x_{\alpha_{\xi}}\mapsto x_{\beta_{\xi}}.\]

Furthermore, define the set \[K_{0}=\{s\in K:|\omega_{g}(s)|\geq 2\},\] where\footnote{Here $B_{\frac{1}{n}}^{K}(s)$ is the ball of radius $\frac{1}{n}$ with center $s$ in $K$, a subset of $\mathds{R}$. } $\omega_{g}(s)=\bigcap_{n\in \omega}\overline{g[B_{\frac{1}{n}}^{K}(s)]}$ is the oscillation of $g$ at $s$, i.e., all the accumulation points of the images (under $g$) of sequences that converges to $s$. Notice that $|\omega_{g}(s)|=1$ if and only if $g$ is continuous at $s$.

Recall that any partial continouos function from $\mathds{R}^{n}$ to $\mathds{R}$ can be extended with a partial function whose domain is a $G_{\delta}$ set. With this, our construction of $\mathds{R}_{\mathcal{S}}$, and the fact that $\mathds{R}_{x_{\alpha_{\xi}}}\neq \mathds{R}_{x_{\beta_{\xi}}}$, we have that the set $K_{0}$ is of size continuum.

Given $s\in K_{0}$ call $a_{s}$, $b_{s}$ two elements in $\omega_{g}(s)$. Without loss of generality, we can assume that $a_{s}<b_{s}$. Let $r\in \mathds{Q}$ such that $a_{s}<r<b_{s}$. Since we only have countably many options, we may shrink $K_0$ in such a way that for all $s\in K_{0}$ the rational number $r$ is the same. Notice that $K_{0}$ still has size continuum. 

Take $t, s\in K_{0}$ such that $t< s$ and take disjoint intervals $I_{t}, I_{s}$ such that $t\in I_{t}$ and $s\in I_{s}$. By the definition of $a_{t}$, $a_{s}$, $b_{t}$ and $b_{s}$ there are $t_{0},t_{1}\in K\cap I_{t}$ and $s_{0}, s_{1}\in I_{s}$ such that $g(t_{0}), g(s_{0})<r<g(t_{1}),g(s_{1})$. Then for the pairs $(t_{0},g(t_{0}))$, $(s_{1}, g(s_{1}))$ we have $t_{0}<s_{1}$ and $g(t_{0})<g(s_{1})$; and for the pair $(t_{1}, g(t_{1}))$, $(s_{0}, g(s_{0}))$ we have $t_{1}<s_{0}$ but $g(t_{1})>g(s_{0})$.

\end{proof}

The above proof works as an example of how to use the technique and paves the way for the following theorem. 

\begin{theorem}
There are $2^{\mathfrak{c}}$ many $\mathfrak{c}$-dense order types inside $\mathds{R}$ that are closed under Turing equivalence.
\end{theorem}

\begin{proof}
Let $\mathcal{S}$ be a collection of Turing degrees such that $\mathds{R}_{\mathcal{S}}$ is layered $\mathfrak{c}$-$2$-entangled.

First, fix two disjoint subsets of size continuum of $\mathcal{S}$, call them $\mathcal{A}$ and $\mathcal{B}$. We have that any one-to-one function between $\mathds{R}_{\mathcal{A}}$ and $\mathds{R}_{\mathcal{B}}$ is going to be a collection of size continuum of disjoint $2$-tuples of $\mathds{R}_{\mathcal{S}}$ such that each entry of the pair has a different Turing degree ($\mathcal{A}$ and $\mathcal{B}$ are disjoint), therefore, there are   $(x_{1}, y_{1})$, $(x_{2}, y_{2})$, $(w_{1}, z_{1})$, $(w_{2}, z_{2})$ such that $x_{1}<x_{2}$ and $y_{1}<y_{2}$, $w_{1}<w_{2}$ but $z_{1}>z_{2}$. In other words, the function cannot be order preserving or order reversing. This shows that $\mathds{R}_{\mathcal{A}}$ and $\mathds{R}_{\mathcal{B}}$ are not order isomorphic.

Now, assume that $\mathcal{A}$ and $\mathcal{B}$ are different subsets of $\mathcal{S}$ with $|\mathcal{A}\Delta \mathcal{B}|=\mathfrak{c}$, where $a\Delta b=(a\setminus b)\cup (b\setminus a)$. Without lost of generality, assume that $\mathds{R}_{\mathcal{A}\setminus \mathcal{B}}$ is of size continuum. Since $\mathds{R}_{\mathcal{A}\setminus \mathcal{B}}$ is disjoint from $\mathds{R}_{\mathcal{B}}$, and $\mathds{R}_{\mathcal{S}}$ is layered $\mathfrak{c}$-$2$-entangled, we have that $\mathds{R}_{\mathcal{A}\setminus \mathcal{B}}$ is not order isomorphic to any subset of $\mathds{R}_{\mathcal{B}}$. Then, there cannot be any order isomorphism between $\mathds{R}_{\mathcal{A}}=\mathds{R}_{\mathcal{A}\cap \mathcal{B}}\cup \mathds{R}_{\mathcal{A}\setminus \mathcal{B}}$ and $\mathds{R}_{\mathcal{B}}$. So, again, this shows that $\mathds{R}_{\mathcal{A}}$ and $\mathds{R}_{\mathcal{B}}$ are not order isomorphic.

Since there are $2^{\mathfrak{c}}$ subsets of size continuum such that $|\mathcal{A}\Delta \mathcal{B}|=\mathfrak{c}$ between any two of them, we have $2^{\mathfrak{c}}$ non-isomorphic order types. 
\end{proof}

It is interesting to wonder which kind of relation have the Turing degrees of the reals in the set created in Theorem \ref{continuum entangled}. There is the possibility that all of them are inside a cone or that all of them have some intricate relationship. Because any Turing degree has, at most, countably many degrees below it, as long as $\neg CH$ holds, there is no way to make these degrees a tower. At the end, the relation depends on how $\mathcal{S}$ was originally. 

Now, it is possible to combine topological or measure arguments with the proof of Theorem \ref{continuum entangled}. For example, assuming that $\max\{cov(\mathcal{M}),cov(\mathcal{N})\}=\mathfrak{c}$, it is possible to produce an antichain directly. Although it is possible for $\mathcal{S}$ to be an antichain to begin with in (in Theorem \ref{continuum entangled}), the following theorem is an interesting application.

\begin{definition}
Denote by $cov(\mathcal{M})$ the least amount of meager sets that are needed to cover $\mathds{R}$. Analogously, $cov(\mathcal{N})$ denotes the least amount of null sets that cover $\mathds{R}$.
\end{definition}

\begin{corollary}
($\max\{cov(\mathcal{M}),cov(\mathcal{N})\}=\mathfrak{c}$) There is a continuum size antichain of Turing degrees $\mathcal{S}$ such that $\mathds{R}_{\mathcal{S}}$ is layered $\mathfrak{c}$-$2$-entangled.
\end{corollary}

\begin{proof}
The proof is analogous to the one above.

Let $\{f_{\alpha}: \alpha<\mathfrak{c}\}$ be an enumeration of all continuous functions from a $G_{\delta}$ subset of $\mathds{R}$ to $\mathds{R}$ and let $\{\mathbf{a}_{\alpha} :\alpha<\mathfrak{c}\}$ an enumeration of all Turing degrees. Given a Turing degree $\mathbf{a}$, we will denote by $\bigvee \mathbf{a}$ all reals that have Turing degree that computes $\mathbf{a}$ or that is computed by $\mathbf{a}$.

Notice that $\bigvee \mathbf{a}$ is both meager and null: it is a well known result that the upper cones of Turing degrees are meager and null  (see \cite{HinmanHier} and \cite{SacksDeg}, respectively), furthermore, the downward cone is countable, so, it is also meager and null.

We will construct the set by recursion.

Assume that we already have $\mathcal{S}_{\alpha}=\{\mathbf{a}_{\beta_{\xi}}: \xi<\alpha\}$. Let \[\beta_{\alpha}=\min\{\beta : \mathds{R}_{\mathbf{a}_{\beta}}\cap(\bigcup_{\xi<\alpha}\bigvee \mathbf{a}_{\beta_{\xi}}\cup B_{\alpha})=\emptyset\}\]
where \[B_{\alpha}=\{f_{\xi}(\overline{v}): \xi<\alpha, \overline{v}\in \bigcup_{\xi<\alpha}\mathds{R}_{\mathbf{a}_{\beta_{\xi}}}\}.\]

Using $\max\{cov(\mathcal{M}),cov(\mathcal{N})\}=\mathfrak{c}$, we know that $\mathds{R}\setminus (\bigcup_{\xi<\alpha} \bigvee \mathbf{a}_{\beta_{\xi}})$ is of size continuum. On the other hand, $\left|B_{\alpha}\right|=|\alpha|<\mathfrak{c}$. This means that $\mathds{R}\setminus (\bigcup_{\xi<\alpha} \bigvee \mathbf{a}_{\beta_{\xi}}\cup B_{\alpha})$ is of size continuum, so there are at most $|\alpha|$  ordinal $\beta$ such that $\mathds{R}_{\mathbf{a}_{\beta}}\cap(\cup_{\xi<\alpha}\mathds{R}_{x_{\xi}}\cup B_{\alpha})\neq\emptyset$.

The set $\mathcal{S}=\mathcal{S}_{\mathfrak{c}}=\{\mathbf{a}_{\beta_{\alpha}}: \alpha<\mathfrak{c}\}$ is the one that we are looking for. By construction, it is clearly an antichain and, using the same proof as above, we have that $\mathds{R}_{\mathcal{S}}$  is layered $\mathfrak{c}$-$2$-entangled, i.e., has the property that for every collection of size continuum of disjoint $2$-tuples of $\mathds{R}_{\mathcal{S}}$ such that each entry of the pair has different degree then there are   $(x_{1}, y_{1})$, $(x_{2}, y_{2})$, $(w_{1}, z_{1})$, $(w_{2}, z_{2})$ such that $x_{1}<x_{2}$ and $y_{1}<y_{2}$, $w_{1}<w_{2}$ but $z_{1}>z_{2}$.

\end{proof}

Another known construction for an entangled set, due to Avraham, is given by taking a set of $\aleph_{1}$ many reals that are Cohen generic with respect a countable model of $H(\theta)$ with $\theta>\mathfrak{c}$, this can be found in \cite{avraham1981martin}. A result of Slaman and Woodin \cite{SlamanManuscript}, proves that the degrees of $5$-generic reals form an automorphism base for the Turing degrees. In Section \ref{Tower Section} we will use an analogous construction to the one done by Avraham to show results related to the existence of an automorphism of the Turing degrees.

\section{Absolute decomposition of functions}\label{Decomposition Section}

In this section the letters $r,x,y$ will represent real numbers (either elements of $\mathds{R}$, $2^{\omega}$ or $\omega^{\omega}$) and the letters $s, l, m, n$ will represent natural numbers (unless otherwise stated).

Fix a countable-to-one function $\varphi:D\subseteq \mathds{R}\rightarrow \mathds{R}$ such that $\langle x, y\rangle \in \varphi$ is arithmetical (or hyperarithmetical) with respect to $x$ and $y$. In this section we will analyze the complexity of the sentence:

\begin{center}
\textit{The function $\varphi$ is the union of countably many monotone functions.}
\end{center}

First, we will analyze how to code a monotone non-decreasing function in a real number. To code a monotone non-increasing function is analogous.

Notice that, given $A\subseteq \mathds{R}$ and a non-decreasing function $g: A\rightarrow \mathds{R}$ we can define a non-decreasing function $f:\mathds{R}\rightarrow [-\infty, \infty]$ as $f(x)=\sup_{y\leq x}\hat{g}(y)$ where $\hat{g}(y)=-\infty$ whenever $x\notin A$ and $\hat{g}(y)=g(y)$ otherwise. Since $f\restrict A=g$, we can assume that all non-decreasing functions are from the real numbers to $[-\infty, \infty]$.

Now, given a non-decreasing function $f$, we know that it can only have countably many discontinuities. If we call the discontinuities $d_{n}$ and we enumerate the rational numbers as $q_{n}$ then we have that:

\[f(x)=\sup\left(\{f(d_{n}): d_{n}\leq x\}\cup \{f(q_{n}): q_{n}\leq x\}\right).\]

To show that the above equality is true, given $x$ there are only two options: either $f$ is continuous at $x$ in which case $f(x)=\sup\{f(q_{n}): q_{n}\leq x\}$ (and, even if there is a sequence of discontinuities that converges at $x$, by the continuity at $x$ the value should be the same) or $f$ is not continuous at $x$ so there is some $n$ such that $x=d_{n}$ and, because $f$ is non-decreasing the above equality is realized.

This means that, given the collection $\{\langle d_{n},f(d_{n})\rangle : n\in \omega\rangle\}\cup \{\langle q_{n}, f(q_{n})\rangle : n\in \omega\}$ we can decode $f$. The numbers $d_{n}$, $f(d_{n})$ and $f(q_n)$ are reals but there are only countably many of them, so we can code them up in a single real, say, $r_{f}$.

Notice that, in this situation, $y=f(x)$ if and only if the formula $\langle x, y\rangle \in f$ is the same as the following arithmetic (in $x$, $y$ and $r_{f}$) formula: \[\forall n, m \left(\begin{array}{l}
(q_n\leq x \leftrightarrow  f(q_{n})\leq y) \vee (d_{m}\geq x \leftrightarrow  f(d_{m})\geq y )\end{array}\right).\]

Furthermore, we know that it is possible for a single real number to code countably many real numbers and, as a result of the above analysis, countably many monotone functions. Given a real number $r$, we will call $f^{r}_{n}$ the $n$-th monotone function that $r$ is coding. In the case that $r$ do not satisfy the required conditions to code countably many monotone functions, we let $f^{r}_{n}$ be the constant $0$ for all $n\in \omega$.

This means that the sentence \textit{the function $\varphi$ is the union of countably many monotone functions}
can be represented by the following $\Sigma_{2}^{1}$ formula:

\[\exists r \forall x,y ( \langle x,y\rangle \notin \varphi \vee \exists n \langle x,y\rangle \in f^{r}_n)\]

With this, we can show the following theorem:

\begin{theorem}\label{Decomposition}
(ZFC)For every $w\in\mathds{R} $ and every countable-to-one function $\varphi:D\subseteq \mathds{R}\rightarrow \mathds{R}$ that is $\Sigma^{1}_{1}[w]$, i.e, such that $\langle x,y \rangle \notin \varphi$ is a $\Pi^1_1[w]$ formula, $\varphi$ is contained in a countable union of monotone functions. In particular, arithmetic and computable functions are the countable union of monotone functions.
\end{theorem}

\begin{proof}
Fix $w\in \mathds{R}$ and a countable-to-one function $\varphi:D\subseteq \mathds{R}\rightarrow \mathds{R}$ that is $\Sigma^{1}_{1}[w]$.

We will show that all models of ZFC  that have $\varphi$ satisfy the sentence \begin{center}
\textit{The function $\varphi$ is the union of countably many monotone functions.}
\end{center}

Given a model $M$ of ZFC, we can take $L^{M}[\varphi,w]$ which is a model of ZFC+CH and has $\varphi$. By a result of Avraham-Rudin-Shelah \cite{abraham1985consistency}, starting from a model of CH, there is a ccc forcing such that the resulting model, say $N$, satisfies that given $A,B\subseteq \mathds{R}$ and $f:A\rightarrow B$ countable-to-one function, we have that $f$ is contained in a countable union of monotone functions.

In particular, since $\varphi\in L^{M}[\varphi,w]\subseteq N$ we have that 

\[N\vDash \exists r \forall x,y ( \langle x,y\rangle \notin \varphi \vee \exists n \langle x,y\rangle \in f^{r}_n)\]

Since the above sentence is $\Sigma^{1}_{2}[w]$ (i.e, $\Sigma^{1}_{2}$ with parameters in $L[\varphi, w]$) , using Shoenfield absolutness Theorem, we know that it is absolute. Therefore we have that

\[M\vDash \exists r \forall x,y ( \langle x,y\rangle \notin \varphi \vee \exists n \langle x,y\rangle \in f^{r}_n).\]

\end{proof}

\begin{corollary}
Given $f: D\subseteq \mathds{R}\rightarrow \mathds{R}$ a Borel function, $f$ is contained in a countable union of monotone functions.
\end{corollary}

\begin{proof}
Fix $f$ a Borel function. Using borel codes, we know that there is $x\in \mathds{R}$ such that $f$ is $\Sigma^{1}_{1}[x]$ (actually, $\Delta^{1}_1[x]$). Now we just need to use Theorem \ref{Decomposition}.
\end{proof}

\begin{theorem}\label{Turing is layered entangled}
Let $E$ be a $2$-$\mathfrak{c}$-entangled set. Then $\mathds{R}_{E}$ is  $2$-$\mathfrak{c}$-layered entangled.
\end{theorem}

\begin{proof}

Towards a contradiction, assume that we have a 2-$\mathfrak{c}$-entangled set $E$ and a one-to-one monotone function $f:D\subseteq \mathds{R}_{E}\rightarrow\mathds{R}_{E}$ such that for continuum many $x\in \mathds{R}_E$, $f(x)$ is not Turning equivalent to $x$.

The idea of the proof is to construct a one-to-one (partial) monotone function $g$ from $E$ to $E$ such that $g(r)$ is not equal to $r$ for continuum many $r\in E$.

To do this, we want to first take $r\in E$ with $\mathds{R}_{r}\cap D\neq \emptyset$ to a Turing equivalent $x\in D$ with a monotone function and then send $f(x)$ to a Turing equivalent element of it in $E$ using a (most likely different) monotone function.

Given $r\in E$ and $x\in \mathds{R}_{r}\subseteq \mathds{R}_{E}$ we can use a Turing functional to send $r$ to $x$. Since there are only countably many of them, we know that continuum many $r\in E$ will use the same Turing functional. We would like to use Lemma \ref{Decomposition} to have a monotone function. Nevertheless, the problem is that Turing functionals may not be one-to-one, so we cannot use the Lemma. Fortunately, we can fix this.

Let $r\in E$ such that $\mathds{R}_{r}\cap D\neq \emptyset$. We know that there exist $e,d\in \omega$ such that $\varphi_{e}^{r}\in D$ and $\varphi_{d}^{\varphi_{e}^{r}}=r$. Since there are only countably many $\langle e, d\rangle$ we know that continuum many $r\in E$ use the same pair. Without loss of generality, we can assume that there are $e,d\in \omega$ such that for all $r\in E$, $\varphi_{e}^{r}\in D$ and $\varphi_{d}^{\varphi_{e}^{r}}=r$.

Let $B_{0}=\{r\in \mathds{R}: \varphi_{d}^{\varphi_{e}^{r}}=r\}$. This set is arithmetic, so, the function $\varphi_{e}^{\_}:B_{0}\rightarrow \varphi_{e}^{\_}[B_{0}]$ is one-to-one (with inverse $\varphi_{d}$) arithmetic function  and, therefore, a $\Sigma^{1}_{1}$ function. Notice that, since $B_{0}$ is not computable, it may be the case that $\varphi_{e}^{\_}\restrict_{B_{0}}$ is not computable. Nevertheless, it will be $\Sigma^{1}_{1}$.

By Lemma \ref{Decomposition}, $\varphi\restrict_{B_{0}}$ is contained in the union of countably many monotone functions. One of these functions, call it $h_{0}:C_{0}\rightarrow h_{0}[C_{0}]$, will satisfy that $A_{0}=C_{0}\cap B_{0}\cap E$ is of size continuum and $h_{0}\restrict_{A_{0}}=\varphi_{e}\restrict_{A_{0}}$. Notice that $h_{0}\restrict_{A_{0}}$ is an strictly (one-to-one) monotone function such that $h_{0}[A_{0}]\subseteq D$ and $A_{0}\subseteq E$ is of size continuum.

Doing an analogous argument, we can find an strictly (one-to-one) monotone function such that $h_{1}[A_{1}]\subseteq f[h_{0}[A_{0}]]$ and $A_{1}\subseteq E$ is of size continuum.

Then, the function $g=h_{0} \circ f \circ h_{1}^{-1}: h_{0}^{-1}[f^{-1}[h_{1}[A_{1}]]]\rightarrow A_{0}$ is a one-to-one monotone function from a subset of $E$ to $E$. Furthermore, we have that, by construction, $deg(g(r))\neq deg(r)$ which means that $g(r)\neq r$.

There is $g'\subseteq g$ such that $g'$ is a one-to-one function between two disjoint size continuum sets of a $\mathfrak{c}$-entangled set. Notice that the set of $g'$ is a continuum size set of disjoint $2$-tuples that contradict the property of $E$ being entangled.
\end{proof}

\begin{corollary}
There is $\mathcal{F}$, a countable collection of monotone functions, such that $\mathcal{F}_x=\mathds{R}_{x}$ for all $x\in \mathds{R}$.
\end{corollary}

It is important to remark that Theorem \ref{Decomposition} is a little unexpected since you can construct computable functions that are not the countable union of measurable monotone functions. The following example is due to Joe Miller and is publish here with his permission.

\begin{lemma}\label{Decomposition counterexample}
(ZFC) There is a computable function that is not the countable union of countably many measurable monotone functions.
\end{lemma}

\begin{proof}

Let $f: D\subseteq (0,1) \rightarrow \mathds{R}$ be such that if $x=(x_{0}, x_{1}, ...)$ has a unique binary expansion given by $(x_{0}, x_{1}, ...)$ then $f(x)=(y_{0}, y_{1},...)$ where $y_{2n+1}=x_{2n}$ and $y_{2n}=x_{2n+1}$. This is a computable function.

Assume that $f\subseteq \bigcup_{n\in \omega} g_{n}$ where $g_{n}$ is a monotone function for all $n$. Since $f=\bigcup_{n\in \omega} g_{n}\cap f$, we will show that there is $n$ such that $g_{n}\cap f$ is nonmeasurable.

Fix $n$ and let $x\in dom(g_{n}\cap f)$, $m\in \omega$ and $\sigma\in 2^{<\omega}$ such that $m+1$ is even,  $|\sigma|=m+1$ and $\sigma$ is an initial segment of the binary representation of $x$. Notice that, for any $y$ such that $\sigma$ is an initial segment of the binary representation of $y$ we have that $|x-y|<2^{-m}$.

For any $y_{00}, y_{01}, y_{10}, y_{11}$ such that $\sigma 00, \sigma 01, \sigma 10, \sigma 11$ are initial segments of their binary representation we have that $y_{00}< y_{01}< y_{10}< y_{11}$ but $f(y_{00})<f( y_{10})< f(y_{01})< f(y_{11})$. This means that for $\tau$ equal to $\sigma 00, \sigma 01, \sigma 10$ or $ \sigma 11$ all the extensions of $\tau$ are not in $dom(f\cap g_{n})$. We can exemplify this as follows: assume that there is $x_{01}\in dom(f\cap g_{n})$ such that $\sigma 01$ is an initial segments of it and that $g_{n}$ is monotone nondecreasing. Then, there is no $y_{10}\in dom (f\cap g)$ such that its binary representation extends $\sigma 10$ since $x_{01} < y_{10} $ but $f(y_{10})<f(x_{01})$. Notice that the size of the interval of all real numbers such that their binary expansion extends $\sigma 10$ has length, at least, $2^{-|\sigma 10|-2}=2^{-m-5}$.

Now, in any of the cases, if we use $\lambda $ to denote the Lebesgue measure and we assume that $dom(f\cap g_{n})$ is measurable, we have that
\[\frac{\lambda ((x-2^{-m},x+ 2^{-m})\cap dom (f\cap g_{n}))}{\lambda ((x-2^{-m}, x+2^{-m}))}\leq \frac{2^{-m+1}-2^{-m-5}}{2^{-m+1}}=1-2^{-6}, \] as long as $m$ is big enough for $(x-2^{-m}, x+2^{-m})\subseteq (0,1)$.

Since the Lebesgue density of any point in $dom(f\cap g_{n})$ cannot be $1$ (it is less than $1-2^{-6} $), we have that $dom(f\cap g_{n})$ is a null set if it is measurable (see Theorem \ref{Lebesgue positive} to read the statement of the Lebesgue density Theorem).

Since $\lambda(dom(f))=1$, it is not a countable union of null sets but, since $dom(f)=\bigcup_{n\in \omega}dom(f\cap g_{n})$ then there is $n$ such that $dom(f\cap g_{n})$ is not measurable.

\end{proof}

Recall that there is a model of ZF with a certain amount of choice where all subsets of reals are measurable. This fact along with Lemma \ref{Decomposition counterexample} seems to contradict Theorem \ref{Decomposition}. Nevertheless, we want to point at the fact that the proof of the Lebesgue Density Theorem seems to use enough choice to create nonmeasurable sets.

\section{Towers of Models}\label{Tower Section}

Following a result in \cite{avraham1981martin}, we can construct an entangled set using a tower of models.

For the convenience of the reader, we reproduce the proof here with more details in some of the steps that are going to be key for us.

\begin{theorem}[Avraham-Shelah]\label{A-S entangled}
Assuming CH there is an uncountable $2$-entangled set of reals.
\end{theorem}

\begin{proof}
We will work with reals as elements of $2^{\omega}$.

Let $\langle r_{\alpha} : \alpha<\omega_{1}\rangle$ be an enumeration of all real numbers.

Let $M_{0}$ be a countable elementary submodel of $H(\aleph_{2})$ (the sets whose transitive closure has size less than $\aleph_{2}$). 

Assuming that we have define $\langle e_{\xi} : \xi<\alpha\rangle$ and $\langle M_{\xi} : \xi\leq\alpha\rangle$, we define $e_{\alpha}$ to be a Cohen generic with respect to $M_{\alpha}$ (which, by the definition of $M_{\xi}$, will be also an elementary submodel of $H(\aleph_{2})$). Now, if $\alpha$ is a limit ordinal, let $M_{\alpha}=\bigcup_{\xi<\alpha} M_{\xi}$.

If $\alpha=\beta+1$, let $M_{\alpha}$ be an elementary extension of $M_{\beta}$ such that $M_{\alpha}$ is also an elementary submodel of $H(\aleph_{2})$, $M_{\beta}\in M_{\alpha}$ looks countable for $M_{\alpha}$ and, if $e_{\beta}=r_{\gamma}$, then $r_{\delta}\in M_{\alpha}$ for all $\delta\leq \gamma$.

The set $E=\left\{e_{\alpha} : \alpha<\omega_{1}\right\}$ is the set that we are looking for. We will see this by contradiction.

Assume that there is an uncountable collection of disjoint pairs \[A=\{(e_{\alpha_{i}}, e_{\beta_{i}}): i<\omega_{1} \}\subseteq E^{2}\] such that whenever $e_{\alpha_{i}}< e_{\alpha_{j}}$ then $e_{\beta_{i}}< e_{\beta_{j}}$.

Without loss of generality, we can assume that given $i<j<\omega_{1}$ we have $\alpha_{i}<\beta_{i}<\alpha_{j}<\beta_{j}$.

Let $N$ be a countable elementary submodel of $H(\aleph_{2})$ such that $A\in N$. Now, let $\gamma_{0}=\omega_{1}\cap N$ and let $\xi_{1}$ be the first ordinal such that $\langle (e_{\alpha_{i}}, e_{\beta_{i}}) : i<\gamma_{0}\rangle\in M_{\xi_{1}}$. This ordinal exists because every countable sequence can be coded by a single real number, therefore, the real and the sequence are added at some stage. By our construction of $E$, for any $\alpha_{i}>\xi_{1}$ the pair $(e_{\alpha_{i}}, e_{\beta_{i}})$ is $M_{\xi_{1}}$ generic for the product of two Cohen forcing.

At this moment we need some notation. The expression $(f,g)\leq (r,s)$ for $f,g\in 2^{<\omega}$ and $r,s\in 2^{\omega}$ means that $f$ is an initial segment of $r$ and $g$ is an initial segment of $s$. Also, in the same context, the expression $T((f,g), (r,s))$ means that for any real number $\overline{f}$ and $\overline{g}$ that extend $f$ and $g$, respectively, we have that $\overline{f}< r$ if and only if $\overline{g}> s$. With this, notice that the set \[\{(f,g): f,g\in 2^{<\omega}\left(\forall i<\gamma_{0} (f,g)\not\leq (e_{\alpha_{i}}, e_{\beta_{i}})\right)\vee \left(\exists i<\gamma_{0} (T((f,g), (e_{\alpha_{i}}, e_{\beta_{i}})))\right) \},\]
which we will call $D$, is dense in $2^{<\omega}\times 2^{<\omega}$ and belongs to $M_{\xi_{1}}$. By genericity of $(e_{\alpha_{j}},e_{\beta_{j}})$ over $M_{\xi_{1}}$, there is $(f,g)\leq (e_{\alpha_{j}},e_{\beta_{j}})$ such that $(f,g)\in D$.

Notice that, since $N$ is an elementary submodel of $H(\aleph_{2})$ and $N\cap A=\langle (e_{\alpha_{i}}, e_{\beta_{i}}) : i<\gamma_{0}\rangle$ we have that for any $j<\omega_{1}$ and for any finite initial segment $(f,g)\leq (e_{\alpha_{j}},e_{\beta_{j}})$ there is $i<\gamma_{0}$ such that $(f,g)\leq (e_{\alpha_{i}},e_{\beta_{i}})$. This implies that there is $\alpha_{j}> \xi_{1}$ and $i<\gamma_{0}$ such that $e_{\alpha_{i}}<e_{\alpha_{j}}$ but $e_{\beta_{i}}>e_{\beta_{j}}$. 

\end{proof}

There is a couple of things that are important to notice.

First of all, it is not necessary to add every single real number to the tower of elementary models. As long as we add all countable subsequences of the generic sequence, we can run the argument.

Also, the dense set $D$ can be described in a $\Sigma^{0}_{2}$ way (actually, $\Delta^{0}_{2}$)  given that you have the sequence $\langle (e_{\alpha_{i}}, e_{\beta_{i}}) : i<\gamma_{0}\rangle$. Therefore, in the face of it, the full genericity of the Cohen generic is not necessary, we can use a $2$-generic real (or a weak $2$-generic, since the set is open dense). Nevertheless, we can do a little better.

\begin{definition}
Given $g\in 2^{\omega}$ we say that $g$ is \emph{(weak) $1$-generic} if and only if for every $\Sigma^{0}_{1}$ (dense) $S\subseteq 2^{<\omega}$ there exists $\sigma\in 2^{<\omega}$,  $\sigma\leq g$, such that either:
\begin{enumerate}
    \item $\sigma\in S$
    \item For all $\tau$ extending $\sigma$, $\tau\notin S$.
\end{enumerate}

\end{definition}

\begin{lemma}\label{Generic amount A-S}
Given the construction in Theorem \ref{A-S entangled}, it is only necessary that $(e_{\alpha_{j}},e_{\beta_{j}})$ is $1$-generic over $\langle (e_{\alpha_{i}}, e_{\beta_{i}}): i<\alpha_{0} \rangle$ for the proof to work.
\end{lemma}

\begin{proof}
Notice that the following subset of $D$ is $\Sigma^{0}_{1}$ over $\langle (e_{\alpha_{i}}, e_{\beta_{i}}): i<\gamma_{0}\rangle$:

\[S=\{(f,g): f,g\in 2^{<\omega}\exists i<\gamma_{0} (T((f,g), (e_{\alpha_{i}}, e_{\beta_{i}}))) \}.\]

Since for every $\alpha_{j}>\gamma_{0}$ we have that $(e_{\alpha_{j}}, e_{\beta_{j}})$ is $1$-generic over $\langle (e_{\alpha_{i}}, e_{\beta_{i}}): i<\alpha_{0} \rangle$. We know that there is $(f,g)\leq (e_{\alpha_{j}}, e_{\beta_{j}})$ such that either $(f,g)\in S$ or for all $(f',g')\geq (f,g)$, $(f',g')\notin S$.

Nevertheless, notice that, by the definition of $\gamma_{0}$, for all $(f,g)\leq (e_{\alpha_{j}}, e_{\beta_{j}})$ there is $i_{(f,g)}< \gamma_{0} $ such that $(f,g)\leq (e_{\alpha_{i_{(f,g)}}}, e_{\beta_{i_{(f,g)}}})$. Since $e_{\alpha_{i_{(f,g)}}}, e_{\beta_{i_{(f,g)}}}$ are not eventually $0$ or $1$ (they are $1$-generic), for each $(f,g)\leq (e_{\alpha_{j}}, e_{\beta_{j}})$ there is $(f', g')\geq (f,g)$ such that $e_{\alpha_{i_{(f,g)}}}\restrict |f'|+1=f'^{\frown}1 $ and $e_{\beta_{i_{(f,g)}}}\restrict |g'|+1= g'^{\frown}0$, then $T((f'^{\frown}0, g'^{\frown}1), (e_{\alpha_{i_{(f,g)}}},e_{\beta_{i_{(f,g)}}}) )$ and $(f'^{\frown}0, g'^{\frown}1)\in S$. This shows that option two of the definition of $1$-generic cannot be satisfied, so there must be $(f,g)\leq (e_{\alpha_{j}}, e_{\beta_{j}})$ such that $(f,g)\in S$.

Therefore, there is $i<\gamma_{0}$ such that $T((e_{\alpha_{j}}, e_{\beta_{j}}),(e_{\alpha_{i}}, e_{\beta_{i}}))$.

\end{proof}

Now, we would like to make an analogous construction of the above proof to show the following theorem:

\begin{theorem}\label{Automorphism thm}
Assuming CH, given any non-trivial automorphism of the Turing degrees, call it $a$, there is no family $G$ of $1$-generics degrees such that:
    \begin{enumerate}
        \item For every degree $\mathbf{y}$ there is a $1$-generic degree over $\mathbf{y}$ in $G$.
        \item For all $\mathbf{x}\in G$, if $\mathbf{x}$ is $1$-generic over $\mathbf{y}$, then there are $x,z\in\mathds{R}$ with $deg(x)=\mathbf{x}$ and $deg(z)=a(\mathbf{x})$ with $(x,z)$ $1$-generic over $\mathbf{y}$.
    \end{enumerate} 
\end{theorem}

By a result of Slaman and Woodin \cite{slaman1986definability, SlamanManuscript}, all automorphism of the Turing degrees can be expressed in an arithmetical way. Therefore, assuming that the theorem is false, there is an automorphism of the Turing degrees that will send a family of $1$-generic to $1$-generics over them. Our strategy for the proof will be to construct a $2$-entangled set of $1$-generics such that the image of the entangled set under the automorphism is inside the entangled set. Then, since all automorphism can be described in an arithmetical way, using the Decomposition Lemma \ref{Decomposition}, we would generate a monotone function from an entangled set to itself, which is a contradiction.

This approach will run into one problem if we use the construction in \cite{avraham1981martin}:  we could have the image of a $1$-generic be inside the first model where the real appear (since, elementarity implies arithmetic elementarity). We will be able to solve this by being more careful with the amount of ZFC that our models satisfy.

Also, it is important to notice that, if we change point $(2)$ of the theorem to demand full genericity (or $n$-generic for all $n\in \omega$) then the theorem is trivially true. This is, again, because any automorphism can be expressed in an arithmetical way. If you have an automorphism $a$ then there is $k\in \omega$ such that, for all $\textbf{x}$, $a(\textbf{x})\leq_{T} \textbf{x}^{(k)}$. Then, $a(\textbf{x})$ cannot be $k+1$-generic over $\textbf{x}$. Using stronger results, we believe that changing point $(2)$ for $3$-generic is trivially true. Either way, this highlights the importance of Lemma \ref{Generic amount A-S} for Theorem \ref{Automorphism thm}.

We will break this proof into a definition and couple of lemmas.

\begin{definition}
A function $f:A\subseteq\mathds{R}\rightarrow \mathds{R}$ is \emph{powerful} if and only if
\begin{enumerate}
    \item $x\equiv_{T} y$ implies $f(x)\equiv_{T} f(y)$. 
    \item There is a family $G$ of $1$-generics degrees such that
    \begin{enumerate}
        \item For all $x\in \mathds{R}$, there is $\mathbf{y}\in G$ and $y\in \mathbf{y}$ such that $y$ is $1$-generic over $x$.
        \item For all $\mathbf{x}\in G$ there is $x\in \mathbf{x}$ such that $x\in A$.
        \item For all $\mathbf{x}\in G$, if $\mathbf{x}$ is a $1$-generic degree over $y$, then there exist $x_1, x_2\in \mathbf{x}$ and $z\equiv_{T} f(x_1)$ such that $(x_2, z)$ is $1$-generic over $y$.
    \end{enumerate} 
\end{enumerate}
\end{definition}

Notice that (c), specifically the fact that $(x_2,z)$ is $1$-generic, implies that if $deg(x)\in G$ then $x\not\equiv_{T} f(x)$.

\begin{lemma}\label{Powerful Entangled}
(CH) Given a powerful function $g$, there is a $2$-entangled set $E$, made of $1$-generics, such that $g(\mathds{R}_{E})\cap \mathds{R}_{E}$ is uncountable. 
\end{lemma}

\begin{proof}

Let $\langle r_{\alpha} : \alpha<\omega_{1}\rangle$ be an enumeration of all real numbers. We will start with some notation: given $A\subset \mathds{R}$ a countable set, let $O_{A}$ be one of the followings:
\begin{enumerate}
    \item The minimal real (in the enumeration) that belongs to the minimal Turing degree (respect to Turing reduction) such that it can compute all elements of $A$.
    \item If the above one is not possible, then $O_{A}$ will be the minimal real (in the enumeration) such that it can compute all elements of $A$  but does not compute $\mathbf{0}'$ (zero jump).
    \item  If neither of the above ones are possible, $O_{A}$ will be the minimal real (in the enumeration) such that it can compute all elements of $A$.
\end{enumerate}

Now, let $M_{0}=H(\aleph_{0})\cup C$, where $C$ is the set containing all computable subsets of $H(\aleph_{0})$. Here, $H(\aleph_{0})$ is the set of the hereditary finite sets. 

Assuming that we have define $\langle e_{\xi} : \xi<\alpha\rangle$ and $\langle M_{\xi} : \xi\leq\alpha\rangle$. If $\alpha=\gamma+2n$ for $\gamma=0$ or $\gamma$ a limit ordinal, we define $e_{\alpha}$ to be the minimal real (with respect the enumeration) such that it is $1$-generic with respect to $O_{\mathds{R}\cap M_{\alpha}}$  and such that $deg(e_{\alpha})\in G$ (the family associated with the powerful function). Now, let $e_{\alpha+1}$  be such that $(e_{\alpha}, e_{\alpha+1} )$ is $1$-generic over $O_{\mathds{R}\cap M_{\alpha}}$ and $e_{\alpha+1}\equiv_{T}g(x_{\alpha})$, where $x_{\alpha}\equiv_{T} e_{\alpha}$. We can do this because $g$ is powerful and $deg(e_{\alpha})\in G$.

Throughout the construction we will do a bookkeeping (as in \cite{solovay1971iterated}) of the countable subsets of $\langle e_{\xi} : \xi<\alpha\rangle$. The objective is that, once we have the sequence $\langle e_{\xi} : \xi<\omega_{1}\rangle$, we can also enumerate all the countable subsets of it as $\langle B_{\xi}: \xi<\omega_{1}\rangle$. Notice that this is possible since every countable subset of $\langle e_{\xi} : \xi<\omega_{1}\rangle$ is a countable subset of $\langle e_{\xi} : \xi<\alpha\rangle$, for some $\alpha<\omega_{1}$ (and by CH).

To finish the construction, if $\alpha$ is a $\delta$ limit ordinal, let $M_{\alpha}=\bigcup_{\xi<\alpha} M_{\xi}\cup C[O_{B_{\delta}}]$ where $ C[O_{B_{\delta}}]$ is the set of all the computable objects from $O_{B_{\delta}}$. Finally, if $\alpha=\beta+1$, let $M_{\alpha}=M_{\beta}\cup C[e_{\beta}]$.

The set $E=\left\{e_{\alpha} : \alpha<\omega_{1}\right\}$ is $2$-entangled. We will see this by contradiction.

Assume that there is an uncountable collection of disjoint pairs \[A=\{(e_{\alpha_{i}}, e_{\beta_{i}}): i<\omega_{1} \}\subseteq E^{2}\] such that whenever $e_{\alpha_{i}}< e_{\alpha_{j}}$ then $e_{\beta_{i}}< e_{\beta_{j}}$.

As in Theorem \ref{A-S entangled}, we can ask that given $i<j<\omega_{1}$ we have $\alpha_{i}<\beta_{i}<\alpha_{j}<\beta_{j}$.

Let $N$ be a countable elementary submodel of $H(\aleph_{2})$ such that $A\in N$. Now, define $\gamma_{0}=\omega_{1}\cap N$ and let $\xi_{1}$ be the first ordinal such that $O_{\Gamma}\in M_{\xi_{1}}$, where $\Gamma=\{ (e_{\alpha_{i}}, e_{\beta_{i}}) : i<\gamma_{0}\}$. By our construction of $E$, for any $\alpha_{i}>\xi_{1}$ the real $(e_{\alpha_{i}}, e_{\beta_{i}})$ is $1$-generic over $O_{\Gamma}$.

Using the same notiation as in Theorem \ref{A-S entangled}, we know that the set \[\{(f,g): f,g\in 2^{<\omega}\left(\forall i<\gamma_{0} (f,g)\not\leq (e_{\alpha_{i}}, e_{\beta_{i}})\right)\vee \left(\exists i<\gamma_{0} (T((f,g), (e_{\alpha_{i}}, e_{\beta_{i}})))\right) \},\]
which we will call $D$, is open dense in $2^{<\omega}\times 2^{<\omega}$.

By Lemma \ref{Generic amount A-S}, we know that there is $\alpha_{i}> \xi_{1}$ and $j<\alpha_{0}$ such that $e_{\alpha_{i}}<e_{\alpha_{j}}$ but $e_{\beta_{i}}>e_{\beta_{j}}$.

This shows that $E$ is $2$-entangled and, by construction, $\mathds{R}_{E}\cap g(\mathds{R}_{E})$ is uncountable.
\end{proof}

\begin{lemma}\label{no powerful func}
(CH) There is no countable-to-one $\Sigma^{1}_{1}$ powerful function.
\end{lemma}

\begin{proof}
Suppose that there exist such a function. Then, by Lemma \ref{Powerful Entangled}, there is an entangled set $E$, made of $1$-generics, such that $g(\mathds{R}_{E})\cap \mathds{R}_{E}$ is uncountable. By Theorem \ref{Turing is layered entangled}, $\mathds{R}_{E}$ is $\sigma$-$2$-entagled.

Since $g$ is a countable-to-one $\Sigma^{1}_{1}$ function, by lemma \ref{Decomposition}, it contains an uncountable monotone function, $g_{1}$, such that $g_{1}(\mathds{R}_{E})\cap \mathds{R}_{E}$ is uncountable.

Nevertheless, by the obsevation made after the definition of powerful, there are uncountably many $1$-generics in $E$ such that $g(x)\not\equiv_{T} x$.

Then we have a strictly monotone function changing uncountably many degrees from an uncountable subset of a $\sigma$-$2$-entangled set to itself. This cannot happen.

\end{proof}

With the above lemmas, we are ready to proof Theorem \ref{Automorphism thm}.

\begin{proof}
By a result in \cite{SlamanManuscript}, every automorphism has an arithmetic representation, this is, a function $f:A\subseteq\mathds{R} \rightarrow \mathds{R}$ such that $f$ is arithmetic, $A$ contains an element from each Turing degree, $x\equiv_{T} y$ if and only if $f(x)\equiv_{T}f(y)$ and the change of degree is done exactly as the automorphism.

Working towards a contradiction, assume that there is an automorphism $a$ and a family $G$ of $1$-generics with the properties enlist in the theorem. The representation $f$ of $a$  is an arithmetic powerful function. By lemma \ref{Powerful Entangled}, this function does not exist.

\end{proof}

\section{Modifications and other possible results}\label{Modification Section}

It is important to make a couple of remarks about the construction of Lemma \ref{Powerful Entangled}. First of all, the construction can be stated purely using Turing degrees. The usage of models was to make the analogy with Theorem \ref{A-S entangled} more direct. Also, we just need the degree of $(e_{\alpha}, e_{\beta})$ to be $1$-generic over the degree of each element of $\mathds{R}\cap M_{\alpha}$. This means that the usage of $O_{\mathds{R}\cap M_{\alpha}}$ was not necessary, nevertheless, it reduces writing. Specifically, because having a degree being $1$-generic over countably many degrees does not mean that you have a single element that is $1$-generic over all of them at the same time (or not in a trivial way), see question \ref{generic same time}.

Either way, it is possible to show that the set that you get is $2$-entangled only assuming that the degree of $(e_{\alpha}, e_{\beta})$ to be $1$-generic over the degree of each element of $\mathds{R}\cap M_{\alpha}$. We could define $e_{\alpha}=x_{\alpha}$ and $e_{\alpha+1}=f(x_{\alpha+1})$. Then, using the $(e,d)$ trick in Theorem \ref{Turing is layered entangled}, we get a monotone function from pairs of $1$-generic degrees over $O_{\Gamma}$ (that is in some $\mathds{R}\cap M_{\alpha}$) to $(e_{\alpha_{i}}, e_{\beta_{i}})$ . The last thing to do will be to change the definition of $D$ depending on whether the function is decreasing or increasing. This description is missing details but, we hope that the reader understands why we choose to make the proof differently.

In a different (but related) topic, do not believe that the strange definition of $O_{A}$ was completely for free. It is important to remark that adding $O_{B_{\alpha}}$, contrary to using $O_{\mathds{R}\cap M_{\alpha}}$, is a key point in the proof since we need the subset $S$ of $D$ of Lemma \ref{Generic amount A-S} to be an open  $\Sigma_{1}^{0}$ set of a single degree.

Finally, if it is the case that $g(x)$ is $1$-generic over $x$ whenever $x$ is $1$-generic, we can show that $0'\notin M_{\alpha}$ for any $\alpha< \omega_{1}$. Making $0'\notin M_{\alpha}$ ensures that for any element in $M_{\alpha}$ all the $1$-generics, respect to that degree, are not in $M_{\alpha}$. Therefore, if we pick $2$-generics instead of $1$-generics, for any $e_{\alpha}$, $g(e_{\alpha})\notin M_{\alpha+1}$. We didn't use this in our proof since $g(e_{\alpha})$ already had all the wanted properties because $g$ is powerful. Nevertheless, this could be handy if the definition of powerful is weakened or for other applications, as the following.

\begin{corollary}
Assuming CH, one of the following is true:

\begin{enumerate}
    \item The relation ``$\textbf{x}$ is $1$-generic over $\textbf{y}$'' is not definable in the Turing degrees.
    \item For any automorphism of the Turing degrees, call it $a$, there is no family $G$ of $1$-generics degrees such that:
    \begin{enumerate}
        \item For every degree $\mathbf{y}$ there is a $1$-generic over $\mathbf{y}$ in $G$.
        \item For all $\mathbf{x}\in G$, if $\mathbf{x}$ and $a(\mathbf{x})$ are $1$-generic degrees over $y$, then there are $x,z\in\mathds{R}$ with $deg(x)=\mathbf{x}$ and $deg(z)=a(\mathbf{x})$ with $(x,z)$ $1$-generic over $y$.
    \end{enumerate} 
\end{enumerate}

\end{corollary}

\begin{proof}
We can modify the definition of $O_{A}$ so that $O_{A}$ computes $0''$. This way, we have that, for any automorphism $a$, $a(deg(O_{A}))=deg(O_{A})$ (see \cite{SlamanManuscript}). Let $f$ be a representation of $a$.

Therefore, if the relation ``$\textbf{x}$ is $1$-generic over $\textbf{y}$'' is definable, we have that if $x$ is $1$-generic over $O_{A}$ then $f(x)$ is generic over $f(O_{A})\equiv_{T}O_{A}$. So, if $1$ and $2$ fails, there is an automorphism $a$ such that $(x, f(x))$ is $1$-generic over $O_{A}$ if $x\in G$.

The rest of the proof follows.
\end{proof}

\begin{corollary}\label{G all 1 gen}
Assuming CH, one of the following is true:

\begin{enumerate}
    \item The relation ``$x$ is $1$-generic over $y$'' is not definable in the Turing degrees.
    \item For any automorphism of the Turing degrees, call it $a$, it is not the case that for all $1$-generic degrees $\mathbf{x}$, if $\mathbf{x}$ and $a(\mathbf{x}) $ are $1$-generic over $y$, then there is $z\in a(\mathbf{x})$ such that $z$ is $1$-generic over $y$ and over $\mathbf{x}$.
\end{enumerate}

\end{corollary}

\begin{proof}
This follows from the previous lemma using the family $G$ to be all $1$-generic degrees.
\end{proof}

Notice that, in both proofs, we are not using that ``y is $1$-generic over x'' is definable but rather that there is a definable relation $R$ such that $xRy$ implies ``$y$ is $1$-generic over $x$'' and such that for all $x$ there is a $y$ with $xRy$. This gives us the following corollary:

\begin{corollary}
Assuming CH, at least, one of the following is true:

\begin{enumerate}
    \item There is no relation $R$ definable over the Turing degrees such that $xRy$ implies ``$y$ is $1$-generic over $x$'' and such that for all $x$ there is a $y$ with $xRy$.
        \item For any automorphism of the Turing degrees, call it $a$, there is no family $G$ of $1$-generics degrees such that:
    \begin{enumerate}
        \item For every degree $\mathbf{y}$ there is a $1$-generic over $\mathbf{y}$ in $G$, say $\mathbf{g}$, such that $\mathbf{y}R\mathbf{g}$.
        \item For all $\mathbf{x}\in G$, if $yR\mathbf{x}$ and $yRa(\mathbf{x})$, for some real $y\in \mathds{R}$, then there are $x,z\in\mathds{R}$ with $deg(x)=\mathbf{x}$ and $deg(z)=a(\mathbf{x})$ with $(x,z)$ $1$-generic over $y$.
    \end{enumerate}
\end{enumerate}
\end{corollary}

Now, a couple of words dedicated to absoluteness are in place. From the proofs it is clear that the use of CH, or of the existance of a tower of models with that property, is essential. Nevertheless, in \cite{SlamanManuscript}, it was shown that the statement ``There is a non-trivial automorphism of the Turing degrees'' was absolute. This does not imply that the above Theorem, or Corollaries, are absolute since there is the possibility that the definition of the $1$-generic relation is not absolute. Indeed, that will depend on the proof (or disproof) of the existence of such a relation.

To see other situations where towers of models can be generated, see the COMA, defined by Hart and Kunen in \cite{KunenCOMA}.

\section{Questions and conclusions}\label{Question section}

In this paper, we studied how do sets closed under Turing equivalance look inside $\mathds{R}$. Definitively, not all the questions were solved. Questioning how these sets behave from an algebraic perspective left some open questions:

\begin{question}
In Theorem \ref{Josiah Subfield}, is the use of $\mathbf{0}'$ necessary? In other words, are there reals $x,y$ such that $\mathbf{0}'\not\leq_{T} deg(x\oplus y)$ and the minimal subfield of $\mathds{R}$ that contains $\mathds{R}_{x}\cup \mathds{R}_{y}$ is not $ \mathds{R}_{\wedge x\oplus y}$? 
\end{question}

\begin{question}
Given a field $F$, is $\mathds{R}_{F}$ also a field?
\end{question}

And, of course, we have questions related to measure like:

\begin{question}
For which subsets $S\subseteq \mathds{R}$ is $\mathds{R}_{S}$ measurable?
\end{question}

Even in the case of order type, which was studied in a much more lengthy way, the advance was not fundamental. So we can ask:

\begin{question}
Which other $\mathfrak{c}$-dense order can be obtained (or cannot be obtained) by a set of the form $\mathds{R}_{S}$?
\end{question}

\begin{question}
Is there a topological or model theoretic characterization of all the order types of the form $\mathds{R}_{S}$?
\end{question}

\begin{question}
Is there a collection of countably many monotone functions, $\mathcal{F}$, such that $\mathcal{F}_x=\mathbf{R}_{\wedge deg(x)}$ for all $x$?
\end{question}

But this approach to subsets of reals showed that it can interact with other important questions in the area of Computability Theory or, as it has been commented to the author, in Descriptive Set Theory.

Our approach to the automorphism problem gave some restrictions to the way that automorphisms interacts with $1$-generics, under the Continuum Hypothesis. The result, as showed in Section \ref{Modification Section}, can be written in multiple ways, but the interactions between the degrees and the reals inside of them are more subtle than they appear.

For example, one way to improve our result, specifically Corollary \ref{G all 1 gen}, will be to show that there is no automorphism such that $a(\mathbf{x})$ is $1$-generic over $\mathbf{x}$ for a big family of $1$-generics. Nevertheless, if we attempt to prove it, the following question arises:

\begin{question}\label{generic same time}
Given a Turing degree $\mathbf{g}$ that is $1$-generic over $\mathbf{a}$ and $\mathbf{b}$, under what conditions is there a real number $z$, with Turing degree $\mathbf{g}$, such that $z$ is $1$-generic over $\mathbf{a}$ and $\mathbf{b}$?
\end{question}

\begin{question}\label{generic automorphism}
Given a sufficiently generic degree $\mathbf{g}$ over $y$, is it true that the image of $\mathbf{g}$ under any nontrivial automorphism of the Turing degrees can compute a  sufficiently generic degree over $y$? Can this degree be $1$-generic over $\mathbf{g}$? 
\end{question}

Finally, Turing reduction is not the only way to classify the real numbers in degrees that show you how much information they bear. For example, we could use the enumeration reduction (here $A\leq_{e} B$ if and only if $A$ is c.e. over $B$) or the constructible reduction (here $x\leq_{c} y$ if and only if $x\in L[y]$, this also have countable degrees if you assume large cardinals).

\begin{question}
Which of the results in this papers are still true, or false, if we change Turing equivalence to enumeration, arithmetic, hyperearithmetic or constructibility equivalence?
\end{question}

In this same line of thought, we can also wonder about the properties of these sets in different spaces.

\begin{question}
What can we say about subsets (or subspaces) of the Hilbert cube that are closed under the continuous degree equivalence (see Miller \cite{MillerContinuous})?
\end{question}

\newpage

\bibliographystyle{abbrv}
\bibliography{biblio}

\end{document}